\documentclass{article}
\usepackage{amsmath}
\usepackage{multicol,color}
\usepackage{float}
\usepackage{soul}
\usepackage{graphicx}
\usepackage{amssymb}
\usepackage{theorem}
\usepackage{fancyhdr}
\usepackage{tikz}
\usepackage{enumerate}
\usepackage[margin=1.1in]{geometry}
\usepackage{pgfplots}
\usepgflibrary{shapes.geometric}
\pgfplotsset{my style/.append style={axis x line=middle, axis y line= middle, xlabel={$x$}, ylabel={$y$}, axis equal }}
\usetikzlibrary{through,calc}
\usetikzlibrary{calc,intersections}
\newcommand\markangle[6][red]{\begin{scope}
\path[clip](#2)--(#3)--(#4);
\fill[color=#1,fill opacity=0.5,draw=#1,name path=circle]
(#3) circle(#6mm);
\end{scope}
\path[name path=line one](#3)--(#2);
\path[name path=line two](#3)--(#4);
\path[
name intersections={of=line one and circle, by={inter one}},
name intersections={of=line two and circle, by={inter two}}](inter one)--(inter two) coordinate[pos=.5] (middle);
\path[name path=bissectrice](#3)--(barycentric cs:#3=-1,middle=1.2);
\path[name intersections={of=bissectrice and circle, by={middleArc}}](#3)--(middleArc) node[pos=1.3] {#5};}

\newcommand*{\TlGreen}[1][1]{%
  \begin{tikzpicture}[every path/.style={thick,fill=lightgray},scale=#1]
  \draw (0,0)    circle (0.07);
  \draw (0.1,0)  circle (0.07);
  \draw[fill=green] (60:0.1) circle (0.07);
  \end{tikzpicture}%
}

\def\disp{\displaystyle}
\def\ve{\varepsilon}
\def\dd{\delta}

\def\lm{\lambda}
\def\O{\Omega}
\def\Tilde{\widetilde}
\def\tilde{\widetilde}
\def\tx{\widetilde{x}}
\def\tu{\widetilde{u}}
\def\ta{\widetilde{a}}
\def\tz{\widetilde{z}}

\def\oa{\bar a}

\def\ox{\bar{x}}
\def\oy{\bar{y}}
\def\oz{\bar{z}}
\def\ov{\bar{v}}
\def\ou{\bar{u}}

\def\gph{\hbox{}}

\def\gg{\gamma}
\def\dn{\downarrow}

\def\Lto{\Longrightarrow}
\def\tto{\rightrightarrows}
\def\st{\stackrel}

\def\Limsup{\mathop{{\rm Lim}\,{\rm sup}}}

\def\hat{\widehat}

\def\tilde{\widetilde}
\def\Tilde{\widetilde}
\def\Bar{\overline}
\def\ra{\rangle}
\def\la{\langle}
\def\ve{\varepsilon}
\def\B{I\!\!B}
\def\h{\hfill\Box}
\def\R{\mathbb{R}}
\def\N{\mathbb{N}}

\def\co{\mbox{\rm co}\,}
\def\gph{\mbox{\rm gph}\,}
\def\epi{\mbox{\rm epi}\,}

\def\dom{\mbox{\rm dom}\,}

\def\dist{\mbox{\rm dist}}

\def\bd{\mbox{\rm bd}\,}

\def\var{\mbox{\rm var}\,}
\def\dn{\downarrow}
\def\O{\Omega}

\def\ph{\varphi}
\def\emp{\emptyset}
\def\st{\stackrel}
\def\oR{\Bar{\R}}
\def\lm{\lambda}

\def\gg{\gamma}

\def\dd{\delta}

\def\al{\alpha}
\def\vt{\vartheta}

\def\be{\beta}

\def\N{I\!\!N}
\def\th{\theta}
\def\vTh{\vartheta}
\newtheorem{theorem}{Theorem}[section]

\newtheorem{proposition}[theorem]{Proposition}
\newtheorem{definition}[theorem]{Definition}
\theoremstyle{plain}{\theorembodyfont{\rmfamily}
}
\theoremstyle{plain}{\theorembodyfont{\rmfamily}
}
\theoremstyle{plain}{\theorembodyfont{\rmfamily}
}
\theoremstyle{plain}{\theorembodyfont{\rmfamily}
\newtheorem{example}[theorem]{Example}}

\theoremstyle{plain}{\theorembodyfont{\rmfamily}
}

\begin{document}
\begin{center}
\vspace*{0.3in}{\bf OPTIMAL CONTROL OF A NONCONVEX PERTURBED SWEEPING PROCESS}\\[3ex]
TAN H. CAO \footnote{Department of Applied Mathematics and Statistics, State University of New York--Korea, Yeonsu-Gu, Incheon, Republic of Korea (tan.cao@stonybrook.edu). Research of this author was partly supported by the Ministry of Science, ICT and Future Planning, Republic of  Korea under the ICT Consilience Creative Program (IITP-2015-R0346-15-1007) supervised by the Institute for Information \& Communications Technology Promotion.}
and B. S. MORDUKHOVICH\footnote{Department of Mathematics, Wayne State University, Detroit, Michigan, USA (boris@math.wayne.edu). Research of this author was partly supported by the US National Science Foundation under grant DMS-1512846, by the US Air Force Office of Scientific Research under grant \#15RT0462, and by the RUDN University Program 5-100.}
\end{center}
\small{\bf Abstract.} The paper concerns optimal control of discontinuous differential inclusions of the normal cone type governed by a generalized version of the Moreau sweeping process with control functions acting in both nonconvex moving sets and additive perturbations. This is a new class of optimal control problems in comparison with previously considered counterparts where the controlled sweeping sets are described by convex polyhedra. Besides a theoretical interest, a major motivation for our study of such challenging optimal control problems with intrinsic state constraints comes from the application to the crowd motion model in a practically adequate planar setting with nonconvex but prox-regular sweeping sets. Based on a constructive discrete approximation approach and advanced tools of first-order and second-order variational analysis and generalized differentiation, we establish the strong convergence of discrete optimal solutions and derive a complete set of necessary optimality conditions for discrete-time and continuous-time sweeping control systems that are expressed entirely via the problem data.\\[1ex]
\vspace*{-0.2in}

\normalsize
\section{Introduction, Problem Formulation, and Discussions}
\setcounter{equation}{0}\vspace*{-0.1in}
The sweeping process was introduced by Jean-Jacques Moreau in the 1970s (see \cite{mor_frict}) in the form
\begin{eqnarray}\label{sp}
\dot x(t)\in-N\big(x(t);C(t)\big)\;\mbox{ a.e. }\;t\in[0,T],
\end{eqnarray}
where $C(t)$ is a (Lipschitz or absolutely) continuous moving {\em convex} set, and where the normal cone $N$ to it is understood in the sense of convex analysis for which Moreau was one of the creators and major players. The original Moreau's motivation came mainly from applications to elastoplasticity, but it has been well recognized over the years that the sweeping process is important for many other applications to various problems in mechanics, hysteresis systems, traffic equilibria, social and economic modelings, etc.; see, e.g., \cite{CT,Kr,KMM,mv2,mi,smb,ve} and the references therein.\vspace*{-0.05in}

Due to the maximal monotonicity of the normal cone operator in convex analysis, the sweeping system \eqref{sp} is described by a dissipative discontinuous differential inclusion and can formally be related to control theory for dynamical systems governed by differential inclusions of the type $\dot x\in F(x)$, which has been broadly developed in variational analysis and optimal control; see, e.g., the books \cite{clsw,m-book2,v} with their extensive bibliographies. However, the results of the latter theory, obtained under certain Lipschitzian assumptions on $F$, are not applicable to the discontinuous sweeping process \eqref{sp}. Moreover, it is well known that the Cauchy problem for \eqref{sp} admits a unique solution, which excludes any optimization and control of the sweeping process in form \eqref{sp} with a given moving set $C(t)$.\vspace*{-0.02in}

The authors of \cite{chhm1} introduced a {\em control version} of the sweeping process by inserting control actions into the moving set $C(t)$ with considering its {\em polyhedral} evolution
\begin{eqnarray}\label{sw-poly}
C(t):=\big\{x\in\R^n\big|\;\la u_i(t),x\ra\le b_i(t),\;i=1,\ldots,m\big\},\;\|u_i(t)\|=1\;\mbox{ for all }\;t\in[0,T],
\end{eqnarray}
where optimal control functions $u_i(t)$ and $b_i(t)$ ought to be selected in order to minimize some cost functional. Formulated in this way optimization models for the controlled sweeping process in \eqref{sp} and \eqref{sw-poly} can be written as optimal control problems for unbounded discontinuous differential inclusions with pointwise state constraints of inequality and equality types, which have never been considered before in optimal control theory. The discrete approximation approach of variational analysis, which significantly extends the one developed in \cite{m95,m-book2} for Lipschitzian differential inclusions, allowed the authors of \cite{chhm2} to derive an adequate set of necessary optimality conditions for such polyhedral sweeping control problems with detailed illustrations of new phenomena by nontrivial examples.\vspace*{-0.05in}

A {\em perturbed version} of the polyhedral controlled sweeping process was considered in \cite{cm1,cm2} in the form
\begin{equation}\label{e:5}
-\dot x(t)\in N\big(x(t);C(t)\big)+f\big(x(t),a(t)\big)\;\mbox{ a.e. }\;t\in[0,T],\quad x(0)=x_{0}\in C(0),
\end{equation}
with controls $a\colon[0,T]\to\R^d$ acting in perturbations and controls $u\colon[0,T]\to\R^n$ acting in the polyhedral moving set generated by the fixed vectors $x^*_i$ as
\begin{equation}\label{e:6}
C(t):=C+u(t)\;\mbox{ with }\;C:=\big\{x\in\R^n\big|\;\la x^*_i,x\ra\le 0\;\mbox{ for all }\;i=1,\ldots,m\big\}.
\end{equation}
The necessary optimality conditions for the controlled sweeping process governed by \eqref{e:5} and \eqref{e:6}, which were derived in  \cite{cm1,cm2} by using discrete approximations and appropriate tools of generalized differentiation, were then applied in \cite{cm2} to an optimal control problem for the {\em crowd motion model in a corridor} \cite{mv2,ve} admitted a sweeping polyhedral description of type \eqref{e:5}, \eqref{e:6}.\vspace*{-0.05in}

Note also that other types of optimization problems for some versions of the sweeping process were considered in the literature without using control parameterizations of the moving sets. Control functions appeared there either in additive perturbations \cite{aht,ac,cmr,et,ri}, or in associated ordinary differential equations \cite{ao,bk}. Necessary optimality conditions for optimal controls in such controlled sweeping models were derived in \cite{ac,bk} by employing some other methods different from \cite{cm1,cm2,chhm1,chhm2} under certain strong smoothness assumptions on the boundaries of compact uncontrolled sweeping sets. In the more recent paper \cite{tols}, the author addressed relaxation issues for sweeping optimization problems with controls in additive perturbations and uncontrolled convex moving sets that were also included in optimization.\vspace*{-0.05in}

In this paper we study a perturbed sweeping process of type \eqref{e:5}, where controls enter both perturbations and the moving set given now in the {\em nonconvex} (and hence {nonpolyhedral}) form as
\begin{equation}\label{e:mset}
C(t):=C+u(t)=\bigcap^m_{i=1}C_i+u(t)\;\mbox{ with }\;C_i:=\big\{x\in\R^n\big|\;g_i(x)\ge 0\big\}\;\mbox{ for all }\;i=1,\ldots,m
\end{equation}
defined by some convex and ${\cal C}^2$-smooth functions $g_i\colon\R^n\to\R$. Since the set $C(t)$ is nonconvex, an appropriate normal cone notion to $C(t)$ should be specified in \eqref{e:5}. For definiteness we choose the {\em proximal} normal cone construction to describe the nonconvex sweeping process under consideration, but actually all the major normal cone notions agree in our setting due the {\em prox-regularity} of the set $C(t)$ under the assumptions made; see Section~2. Besides being of its own theoretical interest and importance, the controlled sweeping process version from \eqref{e:5} and \eqref{e:mset} arises in applications to optimal control of the {\em planar crowd motion model}, which is more adequate for the practical use in comparison of the polyhedral corridor version treated in \cite{cm2}. In fact, this has been our primary motivation for the developments of this paper. The results of the crowd motion applications will be presented in the separate paper \cite{cm-crowd}.\vspace*{-0.05in}

This paper concerns the problem of minimizing the Bolza-type functional
\begin{equation}\label{e:Bolza}
\mbox{minimize }\;J[x,u,a]:=\ph\big(x(T)\big)+\int^T_0\ell\big(t,x(t),u(t),a(t),\dot{x}(t),\dot{u}(t),\dot{a}(t)\big)\,dt
\end{equation}
over the control functions $u(\cdot)\in W^{1,2}([0,T];\R^n)$ and $a(\cdot)\in W^{1,2}([0,T];\R^d)$ generating the corresponding trajectories $x(\cdot)\in W^{1,2}([0,T];R^n)$ of the sweeping differential inclusion \eqref{e:5} with the controlled moving set \eqref{e:mset}, where the time final time $T>0$ and the initial vector $x_0\in\R^n$ are fixed. The precise assumptions on terminal extended-real-valued cost function $\ph\colon\R^n\to\Bar{\R}:=(-\infty,\infty]$ and the running cost/integrand $\ell\colon[0,T]\times\R^{4n+2d}\to\Bar{\R}$ will be formulated in Section~2. In addition to the above, the {\em optimal control problem} $(P)$ under consideration contains the pointwise constraints on the controls
\begin{equation}\label{e:ec}
r_1\le\|u(t)\|\le r_2\;\mbox{ for all }\;t\in[0,T]
\end{equation}
with the fixed constraint bounds $0<r_1\le r_2$. Note that the two inequality constraints in \eqref{e:ec} collapse to the equality one when $r_1=r_2$. The positivity requitement on $r_1$ in \eqref{e:ec} is  motivated by applications.

It is important to emphasize that due the construction of the (proximal) normal cone in \eqref{e:5} and the moving set structure in \eqref{e:mset}, we implicitly have the the pointwise constraints of the other type
\begin{equation}\label{e:mc}
g_i\big(x(t)-u(t)\big)\ge 0\;\mbox{ for all }\;t\in[0,T]\;\mbox{ and }\;i=1,\ldots,m,
\end{equation}
which should be taken into account in the subsequent derivations.\vspace*{-0.05in}

We pursue a {\em threefold goal} in this paper. The {\em first one} is to develop the method of {\em discrete approximations} to study the nonpolyhedral sweeping system in \eqref{e:5}, \eqref{e:mset}, and \eqref{e:ec} as well as the optimal control problem $(P)$ for it by constructing a well-posed sequence of discrete-time control systems such that any sweeping feasible solution can be strongly approximated (in the $W^{1,2}$-norm) by feasible ones for discrete systems and that optimal solutions to the discrete counterparts of $(P)$ strongly converge to an optimal solution for the original sweeping control problem $(P)$. The {\em second goal} is to justify the {\em existence} of optimal solutions to the discrete problems and to derive {\em optimality conditions} for them by employing advanced tools of {\em variational analysis} and {\em generalized differentiation}, which are appropriate and in fact unavoidable in this framework. The {\em final goal} is to use discrete approximations as a {\em vehicle} to establish {\em necessary optimality conditions} for the {\em given local optimal} (in an appropriate sense) solution to $(P)$ by passing to the limit from those obtained for their discrete counterparts. The achievement of the latter goal is also heavily based on employing appropriate generalized differential techniques and, in particular, on {\em second-order subdifferential computations} in variational analysis.\vspace*{-0.05in}

As mentioned, the discrete approximation approach to deriving necessary optimality conditions has been implemented before for Lipschitzian differential inclusions \cite{m95,m-book2} as well as for various versions of the polyhedral sweeping process \cite{cm1,cm2,chhm2}, where the polyhedrality of the moving set was strongly exploited. The nonpolyhedral case of the controlled sweeping process treated below is significantly more involved in comparison with the previous developments in all the major steps of our approach.\vspace*{-0.02in}

The rest of the paper is organized as follows. In Section~2 we recall the needed definitions from variational analysis and generalized differentiation, formulate and discuss the basic assumptions on the initial data of $(P)$, and present some preliminary results widely used below. Section~3 constructs a sequence of discrete approximations for all the constraints of problem $(P)$ simultaneously, without touching optimality so far, and show that {\em any feasible} solution to $(P)$  can be strongly approximated in the $W^{1,2}$-norm by feasible solutions to the discrete-time inclusions that are piecewise linearly extended to the continuous-time interval $[0,T]$ under fairly general assumptions. This line of the strong approximation is continued in Section~5 for {\em local optimal} (in the designated sense) solutions to $(P)$, while the preceding Section~4 is devoted to justifying the existence of global optimal solutions to $(P)$ as well as the definition of ``intermediate" (between weak and strong, the latter included) local minimizers and their relaxation that can be covered by the developed method of discrete approximations.\vspace*{-0.05in}

In Section~6 we start preparations to deriving necessary optimality conditions first for discrete-time sweeping problems and then for the original one. These preparations, which are important for their own sake, include computations of the second-order constructions of generalized differentiation that play a significant role in the subsequent results. Section~7 presents the derivation of necessary optimality conditions for discrete approximation problems by reducing them to nondynamic models of mathematical programming with nonsmooth and nonconvex data together with the usage of generalized differentiation calculus and the second-order computations given above.\vspace*{-0.05in}

Section~8 is a culmination, which establishes a complete set of necessary optimality conditions for the original sweeping control problem $(P)$ by passing to the limit from those obtained in Section~7 for discrete approximations together with rather involved techniques of variational analysis ensuring the appropriate convergence of adjoint trajectories and the validity of the limiting relationships. In Section~9 we present two examples, which are related to practical modeling while illustrating the scheme of applications of the obtained necessary optimality conditions to determine optimal solutions. More important practical applications of the obtained optimality conditions appear in \cite{cm-crowd}. The final Section~10 contains concluding remarks and discusses some directions of the future research.\vspace*{-0.02in}

Throughout the paper we use standard notation of variational analysis and control theory; see, e.g., \cite{m-book2,rw,v}. Let us mention that $B(x,r)$ stands for the closed ball of the space in question centered at $x$ with radius $r>0$, $\N:=\{1,2,\ldots\}$, and $x\st{\ph}{\to}\ox$ means that $x\to\ox$ with $\ph(x)\to\ph(\ox)$.\vspace*{-0.2in}

\section{Basic Definitions, Assumptions, and Preliminaries}
\setcounter{equation}{0}\vspace*{-0.1in}

First we recall some definitions from variational analysis systematically used in what follows. The framework of this paper is Euclidean and finite-dimensional. We refer the reader to the books \cite{clsw,m-book1,rw} for more details on generalized differentiation and related issues of variational analysis and to the excellent survey by Colombo and Thibault \cite{CT} on prox-regularity and its applications.\vspace*{-0.05in}

Let $\O\subset\R^n$ be an nonempty set that is locally closed around $\ox\in\O$, and let $\dist(x;\O):=\inf_{y\in\O}\|x-y\|$ be the distance between $x\in\R^n$ and $\O$. The {\em Euclidean projector} of $x$ onto $\O$ is
$$
\Pi(x;\O):=\big\{w\in\O\big|\;\|x-w\|=\dist(x;\O)\big\},\quad x\in\R^n.
$$
which is nonempty if $x$ is sufficiently close to $\ox$. The {\em proximal normal cone} to $\O$ at $\ox$ is given by
\begin{eqnarray}\label{pnc}
N^P(\ox;\O):=\big\{v\in\R^n\big|\;\exists\,\al>0\;\mbox{ such that }\;\ox\in\Pi(\ox+\al v;\O)\big\},\quad\ox\in\O,
\end{eqnarray}
with $N^P(\ox;\O):=\emp$ if $\ox\notin\O$. Another  geometric construction of generalized differentiation used below is the (basic/limiting/Mordukhovich) {\em normal cone} to $\O$ at $\ox\in\O$ defined by
\begin{eqnarray}\label{e:Mor-nc}
N(\ox;\O):=\big\{v\in\R^n\big|\;\exists\,x_k\to\ox,\;w_k\in\Pi(x;\O),\;\al_k\ge 0\;\mbox{ s.t. }\;\al_k(x_k-w_k)\to v\;\mbox{ as }\;k\to\infty\big\}
\end{eqnarray}
with $N(\ox;\O):=\emp$ if $\ox\notin\O$. In contrast to the proximal normal cone \eqref{pnc}, the limiting one \eqref{e:Mor-nc} and the corresponding subdifferential and coderivative constructions for nonsmooth functions and set-valued mappings, which are generated by \eqref{e:Mor-nc} and are presented in Section~6, enjoy {\em full calculi} in general settings that are  based on the {\em variational/extremal principles} of variational analysis; see, e.g., \cite{m-book1,rw}. There is the following relationship between the limiting and proximal normal cone notions:
\begin{eqnarray*}
N(\ox;\O)=\big\{v\in\R^n\big|\;\exists\,x_k\to\ox,\;v_k\to v\;\mbox{ with }\;v_k\in N^P(x_k;\O)\;\mbox{ for all }\;k\in\N\big\}.
\end{eqnarray*}
If the set $\O$ is convex, both constructions \eqref{pnc} and \eqref{e:Mor-nc} reduce to the classical normal cone of convex analysis. But the convexity of $\O$ is not the only case when $N^P(\ox;\O)=N(\ox;\O)$, and thus we can combine nice properties of both cones; in particular, the convexity of \eqref{pnc} and the rich calculus for \eqref{e:Mor-nc}.\vspace*{-0.05in}

It has been well realized in variational analysis that the cones \eqref{pnc} and \eqref{e:Mor-nc} agree for a remarkable class of nonconvex sets introduced in variational analysis by Poliquin and Rockafellar \cite{pr} under the name of {\em prox-regularity}. In fact, this notion was first developed by Federer \cite{f} in geometric measure theory under the name of ``sets with positive reach." The reader can find more information in \cite{CT} with its abounded bibliographies therein. Besides many other applications, prox-regular moving sets have been used in the sweeping process theory; see, e.g., \cite{ac,CT,et,t,ve}. Our main attention is paid to {\em uniformly} prox-regular sets, the notion that was probably first developed by Canino \cite{can} in the study of geodesics.\vspace*{-0.1in}

\begin{definition}{\bf (uniform prox-regularity).}\label{Def-2} Let $\O$ be a closed subset of $\R^n$, and let $\eta>0$. Then $\O$ is $\eta$-{\sc-prox-regular} if for all $x\in\bd\O$ and $v\in N^P(x;\O)$ with $\|v\|=1$ we have $B(x+\eta v,\eta)\cap\O=\emp$. Equivalently, the $\eta$-prox-regularity of $\O$ can be defined via the validity of the estimate
$$
\la v,y-x\ra\le\dfrac{\|v\|}{2\eta}\|y-x\|^2\;\mbox{ for all }\;y\in\O,\;x\in\bd\O,\;\mbox{ and }\;v\in N^P(x;\O).
$$
\end{definition}\vspace*{-0.05in}
Recall that any closed convex subset in $\R^n$ is $\infty$-prox-regular and that for every $\eta>0$ the well-defined Euclidean projection operator $\Pi(x;\O)$ is single-valued whenever $\dist(x;\O)<\eta$.\vspace*{-0.02in}

Throughout the paper we impose the following assumptions on the given data of the optimal control problem $(P)$ ensuring, in particular, that for each $t\in[0,T]$ the controlled moving set $C(t)$ in \eqref{e:mset} is uniformly prox-regular, and thus the proximal and limiting normal cones agree for it. This allows us to use the normal cone notation ``$N$" in \eqref{e:5} in the rest of the paper and employ the results available in variational analysis for either one of these cones. Here are our {\em standing assumptions}:\vspace*{-0.05in}

{\bf (H1)} The perturbation mapping $f\colon\R^n\times\R^d\to\R^n$ in \eqref{e:5} is continuous on $\R^n\times\R^d$ and locally Lipschitzian with respect to the first argument, i.e., for every $\ve>0$ there is a constant $K>0$ such that
\begin{equation}\label{e:Lips}
\|f(x,a)-f(y,a)\|\le K\|x-y\|\;\mbox{ whenever }\;(x,y)\in B(0,\ve)\times B(0,\ve),\;a\in\R^d.
\end{equation}
Furthermore, there is a constant $M>0$ ensuring the growth condition
\begin{equation}\label{e:growth-con}
\|f(x,a)\|\le M(1+\|x\|)\;\mbox{ for any }\;x\in\bigcup_{t\in[0,T]}C(t),\;a\in\R^d.
\end{equation}\vspace*{-0.15in}

{\bf (H2)} There exist positive constants $M_1,M_2$, and $M_3>0$ such that the functions $g_i(\cdot)$, $i=1,\ldots,m$, are twice continuously differentiable (${\cal C}^2$-smooth) satisfying the estimates
\begin{equation}\label{e:bound-grad}
M_1\le\|\nabla g_i(x)\|\le M_2\;\mbox{ and }\;\|\nabla^2 g_i(x)\|\le M_3\;\mbox{ for all }\;x\in V_i\mbox{ for all }\;x\in V_i.
\end{equation}\vspace*{-0.15in}

{\bf (H3)} There exist positive numbers $\be$ and $\rho$ such that
\begin{equation}\label{e:w-inverse-triangle-in}
\sum_{i\in I_\rho(x)}\lm_i\|\nabla g_i(x)\|\le\be\Big\|\sum_{i\in I_\rho(x)}\lm_i\nabla g_i(x)\Big\|\;\mbox{ for all }\;x\in C\;\mbox{ and }\;\lm_i\ge 0,
\end{equation}
where the index set for the perturbed constraints is defined by
\begin{equation}\label{e:rho-index}
I_\rho(x):=\big\{i\in\{1,\ldots,m\}\big|\;g_i(x)\le\rho\big\}.
\end{equation}\vspace*{-0.2in}

{\bf (H4)} The terminal cost $\ph\colon\R^n\to\Bar{\R}$ is lower semicontinuous (l.s.c.), while the running cost $\ell$ in \eqref{e:Bolza} is such that $\ell_t:=\ell(t,\cdot)\colon\R^{4n+2d}\to\Bar{\R}$ is l.s.c.\ for a.e.\ $t\in[0,T]$, bounded from below on bounded sets, and $t\mapsto\ell(t,x(t),u(t),a(t),\dot x(t),\dot u(t),\dot a(t))$ is summable on $[0,T]$ for each feasible triple $(x(t),u(t),a(t))$.\vspace*{-0.02in}

Observe that the simultaneous validity of \eqref{e:bound-grad} and \eqref{e:w-inverse-triangle-in} imply the {\em positive linear independence} of the gradients  $\nabla g_i(x)$ of the active inequality constraints on $C$, and that it reduces to the validity of the classical Mangasarian-Fromovitz constraint qualification on $C$ in the setting under consideration.\vspace*{-0.02in}

The following proposition is due to the result of Venel \cite[Proposition~2.9]{ve}.\vspace*{-0.1in}

\begin{proposition}{\bf (uniform prox-regularity of the moving set).}\label{Pr:prox-reg} Under the validity of {\rm(H2)} and {\rm(H3)} we have that for each $t\in[0,T]$ the set $C(t)$ is $\eta$-prox-regular with $\eta=\dfrac{\al}{M_3\be}\cdot$
\end{proposition}\vspace*{-0.1in}
{\bf Proof.} It follows from \cite[Proposition~2.9]{ve} that the set $C$ in \eqref{e:mset} is $\eta$-prox-regular with the modulus $\eta$ defined in the proposition. Thus it holds for the moving set $C(t)=C+u(t)$ as a translation of $C$. $\h$\vspace*{-0.02in}

Proposition~\ref{Pr:prox-reg} allows us to verify the next proposition based on the well-posedness result by Edmond and Thibault taken from \cite[Theorem~1]{et}.\vspace*{-0.1in}

\begin{proposition}{\bf(existence and estimates for sweeping trajectories).}\label{Pr:1} Consider the perturbed sweeping process \eqref{e:5} with the fixed controls let $u(\cdot)\in W^{1,2}([0,T];\R^n)$ and $a(\cdot)\in W^{1,2}([0,T];\R^d)$ under the validity of {\rm(H1)--(H3)} with the constant $M>0$ taken from \eqref{e:growth-con}. There there is the unique solution $x(\cdot)\in W^{1,2}([0,T];\R^n)$ to \eqref{e:5} generated by the controls $(u(\cdot),a(\cdot))$. Furthermore, we have
$$
\|x(t)\|\le l:=\|x_0\|+e^{2MT}\left(2MT(1+\|x_0\|)+\int^T_0\|\dot{u}(s)ds\|\right)\;\mbox{ for all }\;t\in[0,T],
$$
\begin{equation}\label{e:boundedness}
\|\dot{x}(t)\|\le 2(1+l)M+\|\dot{u}(t)\|\;\mbox{ a.e. }\;t\in[0,T].
\end{equation}
\end{proposition}\vspace*{-0.05in}
{\bf Proof.} With the fixed pair $(u(\cdot),a(\cdot))$, both existence and estimate statements of the theorem follows from \cite[Theorem~1]{et} under the validity of (H1), the uniform prox-regularity of $C(t)$, and the property
$$
\big|\dist\big(y;C(t)\big)-\dist\big(y;C(\tau)\big)\big|\le\|v(t)-v(\tau)\|\;\mbox{ for all }\;t,\tau\in[0,T]\;\mbox{ with }\;v(t):=\int^t_0\|\dot{u}(s)\|ds
$$
and with the chosen $W^{1,2}$ control $u(\cdot)$. The latter fact was proved in \cite[Proposition~1]{cm1}, while the uniform prox-regularity of $C(t)$ follows from Proposition~\ref{Pr:prox-reg} under assumptions (H2) and (H3). $\h$ \vspace*{-0.2in}

\section{Discrete Approximation of Nonconvex Sweeping Process}
\setcounter{equation}{0}\vspace*{-0.1in}

This section deals with the constrained nonconvex sweeping process described by \eqref{e:5}, \eqref{e:mset}, and \eqref{e:ec} without considering its optimization. In what follows we construct a sequence of discrete-time counterparts of the constrained sweeping process in such a way that any feasible triple to the continuous-time process (including controls $u(t)$, $a(t)$ and the corresponding trajectory $x(t)$ satisfying rather unrestrictive conditions) can be {\em strongly} approximated in the $W^{1,2}$-norm  by feasible solutions to the discrete-time systems that are piecewise linearly extended on the continuous-time interval $[0,T]$.\vspace*{-0.05in}

The first step of this procedure is used in all the results presented below. Unifying the control and state variables, we introduce the triple $z:=(x,u,a)\in\R^n\times\R^n\times\R^d$ and show that \eqref{e:5} with $C(t)$ from \eqref{e:mset} can be written in the usual form of differential inclusions with respect to the new variable $z$. Indeed, define the set-valued mapping $F\colon\R^n\times\R^n\times\R^d\tto\R^n$ by
\begin{equation}\label{e:vel-map}
F(z)=F(x,u,a):=N(x-u;C)+f(x,a)
\end{equation}
and deduce from \cite[Proposition~2.8]{ve} that $F$ admits the explicit representation
\begin{equation}\label{e:F-rep}
F(z)=\Big\{-\sum_{i\in I(x-u)}\lm_i\nabla g_i(x-u)\Big|\;\lm_i\ge 0\Big\}+f(x,a)
\end{equation}
via the index set of active constraints
\begin{equation}\label{e:a-index}
I(y):=\big\{i\in\{1,\ldots,m\}\big|\;g_i(y)=0\big\}
\end{equation}
at the point $y:=x-u\in C$. Then we can equivalently rewrite the sweeping inclusion \eqref{e:5} with the set $C$ from \eqref{e:mset} in the following equivalent form involving one variable $z\in\R^n\times\R^n\times\R^d$:
\begin{equation}\label{e:diff-incl}
-\dot{z}(t)\in F\big(z(t)\big)\times\R^n\times\R^d\;\mbox{ a.e. }\;t\in[0,T]\;\mbox{ with }\;z(0):=\big(x_0,u(0),a(0)\big),\;x_0-u(0)\in C,
\end{equation}
where the last condition means that $g_i(x_0-u(0))\ge 0$ as $i=1,\ldots,m$. Proposition~\ref{Pr:1} tells us that the Cauchy problem in \eqref{e:diff-incl} has solutions in the class of $W^{1,2}$ functions $z(t)=(x(t),u(t),a(t))$ on $[0,T]$.\vspace*{-0.02in}

For each $k\in\N$ we define the discrete partitions of $[0,T]$ by
\begin{equation}\label{e:DP}
\Delta_k:=\{0=t^k_0<t^k_1<\ldots<t^k_k\}\;\mbox{ with }\;h_k:=t^k_{j+1}-t^k_j\dn 0\;\mbox{ as }\;k\to\infty.
\end{equation}
The next theorem justifies the desired strong discrete approximation for a large class of feasible solutions to the continuous-time inclusions \eqref{e:diff-incl} generated by the nonconvex sweeping set \eqref{e:mset} with perturbing the pointwise control constraints \eqref{e:ec} and keeping the implicit state constraints \eqref{e:mc}.\vspace*{-0.1in}

\begin{theorem}{\bf (strong discrete approximation of feasible solutions to the constrained sweeping process).}\label{Th:DA} Let the assumptions in {\rm(H1)--(H3)} be satisfied. Consider a triple $\oz(\cdot)=(\ox(\cdot),\ou(\cdot),\oa(\cdot))$ satisfying \eqref{e:diff-incl}, \eqref{e:mset}, and \eqref{e:ec} as well as the following properties, which all hold whenever $\oz\in W^{2,\infty}([0,T])$: the sweeping inclusion \eqref{e:diff-incl} is fulfilled for $\oz(\cdot)$ at the partition points $t^k_j$ from \eqref{e:DP} as $j=0,\ldots,k-1$ with the right-side derivative at $t_0=0$ and for some constant $\mu>0$ independent of $k$ we have
\begin{equation}\label{e:estimates}
\begin{aligned}
&\sum^{k-1}_{j=0}\left\|\dfrac{\ox(t^k_{j+1})-\ox(t^k_j)}{h_k}-\dot{\ox}(t^k_j)\right\|\le\mu,\quad\left\|\dfrac{\ou(t^k_1)-\ou(t^k_0)}{h_k}\right\|\le\mu,\\
&\sum^{k-2}_{j=0}\left\|\dfrac{\ou(t^k_{j+2})-\ou^k(t^k_{j+1})}{h_k}-\dfrac{\ou(t^k_{j+1})-\ou(t^k_j)}{h_k}\right\|\le\mu.
\end{aligned}
\end{equation}
Then there exist a sequence of piecewise linear functions $z^k(t):=(x^k(t),u^k(t),a^k(t))$ on $[0,T]$ and a sequence of positive numbers $\ve_k\le2h_k\mu e^K\dn0$ as $k\to\infty$ with the constant $K>0$ taken from \eqref{e:Lips} so that $(x^k(0),u^k(0),a^k(0))=(x_0,\ou(0),\oa(0))$,
\begin{equation}\label{e:u-discrete-con}
r_1-\ve_k\le\|u^k(t^k_j)\|\le r_2+\ve_k,\quad t^k_j\in\Delta_k,
\end{equation}
\begin{equation}\label{e:x-discrete-con}
x^k(t)=x^k(t_j^k)-(t-t_j^k)v^k_j,\;x^k(0)=x_0,\;t^k_j\le t\le t^k_{j+1}\;\;\mbox{with}\;\;v^k_j\in F(z^k(t^k_j))
\end{equation}
for $j=0,\ldots,k-1$ and that $z^k(\cdot)\to\oz(\cdot)$ in the $W^{1,2}$-norm, i.e.,
\begin{equation}\label{e:W-convergence}
z^k(t)\to\oz(t)\quad\mbox{uniformly on}\quad[0,T]\quad\mbox{and}\quad\int^T_0\left\|\dot{z}^k(t)-\dot{\oz}(t)\right\|^2dt\to 0\quad\mbox{as}\quad k\to\infty.
\end{equation}
Moreover, with $\Tilde\mu:=\max\{3\mu(1+4KT)e^K,4\mu(e^K+1)\}$ we have the estimates
\begin{equation}\label{e:BV}
\left\|\dfrac{u^k(t^k_1)-u^k(t^k_0)}{h_k}\right\|\le\Tilde\mu\;\mbox{ and }\;{\rm var}\big(\dot{u}^k;[0,T]\big)\le\Tilde\mu,
\end{equation}
where ``{\rm var}" denotes the total variation on $[0,T]$ of the derivative function in \eqref{e:BV}.
\end{theorem}\vspace*{-0.1in}
{\bf Proof}. It mainly follows the lines in the proof of \cite[Theorem~3.1]{cm1}, although the problems under consideration in \cite{cm1} and here are significantly different. We just present the major constructions in the new setting of this theorem. Define first a sequence of piecewise linear functions $y^k(\cdot):=(y^k_1(\cdot),y^k_2(\cdot),y^k_3(\cdot))$ on $[0,T]$ via their values at the mesh points of $\Delta_k$ by
$$
\big(y^k_1(t^k_j),y^k_2(t^k_j),y^k_3(t^k_j)\big):=\big(\ox(t^k_j),\ou(t^k_j),\oa(t^k_j)\big)\;\mbox{ for all }\;j=0,\ldots,k\;\mbox{and }\;k\in\N.
$$
Define further $w^k(t)=(w^k_1(t),w^k_2(t),w^k_3(t)):=\dot{y}^k(t)$ as piecewise constant and right continuous functions on $[0,T]$ and easily deduce from these constructions that
$$
y^k(\cdot)\to\oz(\cdot)\;\;\mbox{uniformly on}\;\;[0,T]\;\;\mbox{and}\;\;w^k(\cdot)\to\dot{\oz}(\cdot)\;\;\mbox{strongly in}\;\;L^2([0,T];\R^{2n+d}).
$$
For each fixed $k\in\N$, denote $a^k(t):=y^k_3(t)$ on $[0,T]$ and use for simplicity the notation $t_j:=t^k_j$ as $j=1,\ldots,k$ in what follows. We construct the desired trajectories $x^k(t)$ of \eqref{e:x-discrete-con} by induction. Suppose that the value of $x^k(t_j)$ is known and define the vectors
$$
u^k(t_j):=x^k(t_j)-y^k_1(t_j)+y^k_2(t_j)=x^k(t_j)-\ox(t_j)+\ou(t_j),\quad j=0,\ldots,k,
$$
for which we clearly have the relationships
\begin{equation}\label{e:xk-uk}
x^k(t_j)-u^k(t_j)=\ox(t_j)-\ou(t_j)\;\mbox{ as }\;j=0,\ldots,k\;\mbox{ and}
\end{equation}
$$
r_1-\ve_k\le\|u^k(t_j)\|\le r_2+\ve_k\;\mbox{ with }\;\ve_k:=\|x^k(t_j)-\ox(t_j)\|.
$$
This yields $g_i(x^k(t_j)-u^k(t_j))=g_i(\ox(t_j)-\ou(t_j))\ge 0$ for $i=1,\ldots,m$ and thus $x^k(t_j)-u^k(t_j)\in C$, which shows that the values $F(z^k(t_j))=F(x^k(t_j),u^k(t_j),a^k(t_j))$ of the mapping \eqref{e:vel-map} are well defined whenever $j=0,\ldots,k$. Set $x^k(t):=x^k(t_j)-(t-t_j)v^k_j$ with $v^k_j\in\Pi(-w^k_{1j};F(z^k(t_j)))$ for all $t\in[t_j,t_{j+1})$ and $j=0,\ldots,k-1$ and then deduce from \eqref{e:F-rep} and \eqref{e:xk-uk} that
$$
F\big(x^k(t_j),u^k(t_j),a^k(t_j)\big)=F\big(\ox(t_j),\ou(t_j),\oa(t_j)\big)+f\big(x^k(t_j),\oa(t_j)\big)-f\big(\ox(t_j),\oa(t_j)\big).
$$
Employing now the arguments similar to the proof of \cite[Theorem~3.1]{cm1}, we readily verify that the triples $z^k(t)=(x^k(t),u^k(t),a^k(t))$ for $t\in[0,T]$ and $k\in\N$ constructed above satisfy all the conclusions of this theorem and thus complete the proof of this theorem. $\h$\vspace*{-0.2in}

\section{Existence of Optimal Solutions and Local Relaxation}
\setcounter{equation}{0}\vspace*{-0.1in}

This section starts the study of the entire sweeping optimal control problem $(P)$ formulated in Section~1, not only its feasible solutions. First we establish the existence of optimal solutions to $(P)$ with adding the convexity of the running cost $\ell$ of \eqref{e:Bolza} in the velocity variables to our standing assumptions.\vspace*{-0.1in}

\begin{theorem}{\bf (existence of optimal solutions to the sweeping optimal control problem).}\label{Th:existence} Let $($P$)$ be the optimal control problem defined in Section~{\rm 1} considering in the equivalent form \eqref{e:diff-incl} of the sweeping differential inclusion over all the $W^{1,2}([0,T])$ triples $z(\cdot)=(x(\cdot),u(\cdot),a(\cdot))$. In addition to the standing assumptions in {\rm(H1)--(H4)}, suppose that the integrand $\ell$ in \eqref{e:Bolza} is convex with respect to the velocity variables $(\dot{x},\dot{u},\dot{a})$ and that along a minimizing sequence of $z^k(\cdot)=(x^k(\cdot),u^k(\cdot),a^k(\cdot))$ as $k\in\N$ we have that $\ell(t,\cdot)$ is majorized by a summable function, $\{\dot{u}^k(\cdot)\}$ is bounded in $L^2([0,T];\R^n)$ and $\{a^k(\cdot)\}$ is bounded in $W^{1,2}([0,T];\R^d)$. Then problem $(P)$ admits an optimal solution in $W^{1,2}([0,T];\R^{2n+d})$.
\end{theorem}\vspace*{-0.1in}
{\bf Proof.} It follows from Proposition~\ref{Pr:1} that the set of feasible solutions to $(P)$ is nonempty. Take the minimizing sequence $(x^k(\cdot),u^k(\cdot),a^k(\cdot))\in W^{1,2}([0,T];\R^{2n+d})$ in $(P)$ from the formulation of the theorem, where $x^k(\cdot)$ is uniquely generated by $(u^k(\cdot),a^k(\cdot))$ in Proposition~\ref{Pr:1}.
The imposed boundedness assumptions on $\{(u^k(\cdot)$ and $a^k(\cdot))\}$ yield by standard functional analysis that the sequence $\{(\dot u^k(\cdot),\dot a^k(\cdot))\}$ is weakly compact in $L^2([0,T];\R^{n+d})$, and so we have--by passing to subsequences if necessary--the weak convergence $\dot u^k(\cdot)\to\vTh^u(\cdot)$ and $\dot a^k(\cdot)\to\vTh^a(\cdot)$ in $L^2([0,T];\R^n)$ and $L^2([0,T];\R^d)$, respectively, for some functions $\vTh^u(\cdot)$ and $\vTh^a(\cdot)$ from the corresponding spaces. Due to pointwise constraints \eqref{e:ec} and the boundedness of $\{a^k(0)\}$, suppose without loss of generality that $u^k(0)\to u_0$ and $a^k(0)\to a_0$ as $k\to\infty$ for some $u_0\in\R^n$ and $a_0\in\R^d$. Then $(u^k(\cdot),a^k(\cdot))\to(\ou(\cdot),\oa(\cdot))$ in the norm of $W^{1,2}([0,T];\R^{n+d})$ for the functions $\ou(\cdot)\in W^{1,2}([0,T];\R^n)$ and $\oa(\cdot)\in W^{1,2}([0,T];\R^d)$ defined by
\begin{equation}\label{e:16}
\ou(t):=u_0+\int^t_0\vartheta^u(s)ds\;\mbox{ and }\;\oa(t):=a_0+\int^t_0\vartheta^a(s)ds,
\end{equation}
which implies, in particular, that $\ou(\cdot)$ satisfies the constraints in \eqref{e:ec}. Furthermore, it follows from estimate \eqref{e:boundedness} in Proposition~\ref{Pr:1} that the sequence $\{\dot x^k(\cdot)\}$ is bounded in $L^2([0,T];\R^n)$, and hence we get $\vTh^x(\cdot)\in L^2([0,T];\R^n)$ for which $\dot x^k(\cdot)\to\vTh^x(\cdot)$ weakly in $L^2([0,T];\R^n)$ along a subsequence. It yields the convergence $x^k(\cdot)\to\ox(\cdot)$ in the norm topology of $W^{1,2}([0,T];\R^n)$ for $\ox(\cdot)\in W^{1,2}([0,T];\R^n)$ defined by
$$
\ox(t):=x_0+\int^t_0\vTh^x(s)ds,\quad t\in[0,T].
$$\vspace*{-0.15in}

The next step is to verify that the limiting triple $\oz(\cdot)=(\ox(\cdot),\ou(\cdot),\oa(\cdot))$ satisfies the differential inclusion \eqref{e:diff-incl} with mapping $F(z)$ given in \eqref{e:vel-map}. Since the derivative sequences $\{\dot z^k(\cdot)\}$ converges to $\dot\oz(\cdot)$ weakly in $L^2([0,T];\R^{2n+d})$, the classical Mazur theorem ensures the strong convergence to $\dot\oz(\cdot)$ in this space of some sequence of convex combinations of the functions $\dot z^k(t)$. Thus there is a subsequence of these convex combinations that converges to $\dot\oz(t)$ for a.e. $t\in[0,T]$. It follows from the above that there exist a function $\nu:\N\to\N$ and a sequence of real numbers $\{\alpha(k)_j|\;j=k,\ldots,\nu(k)\}$ such that
$$
\alpha(k)_j\ge 0,\;\sum_{j=k}^{\nu(k)}\alpha(k)_j=1,\;\mbox{ and }\;\sum_{j=k}^{\nu(k)}\alpha(k)_j\dot{z}^j(t)\to\dot{\oz}(t)\;\mbox{ a.e. }\;t\in [0,T]
$$
as $k\to\infty$. By taking into account that $\ox(t)-\ou(t)=\disp\lim_{k\to\infty}({x}^{k}(t)-{u}^{k}(t))\in C$ we get
\begin{equation*}
\begin{aligned}
-\dot{\ox}(t)-&f\big(\ox(t),\oa(t)\big)=\lim_{k\to\infty}\Big(-\sum_{j=k}^{\nu(k)}\alpha(k)_j\dot{x}^j(t)-\sum_{j=k}^{\nu(k)}\alpha(k)_j f\big(x^j(t),a^j(t)\big)\Big)\\
=&\lim_{k\to\infty}\Big(-\sum_{j=k}^{\nu(k)}\sum_{i\in I(x^j(t)-u^j(t))}\al(k)_j\lm^j_i\nabla g_i\big(x^j(t)-u^j(t)\big)\Big)\\
=&\lim_{k\to\infty}\Big(-\sum_{j=k}^{\nu(k)}\sum_{i\in I(\ox(t)-\ou(t))}\al(k)_j\lm^j_i\nabla g_i\big(x^j(t)-u^j(t)\big)\Big)\\
=&\lim_{k\to\infty}\Big(-\sum_{i\in I(\ox(t)-\ou(t))}\sum_{j=k}^{\nu(k)}\al(k)_j\lm^j_i\nabla g_i\big(x^j(t)-u^j(t)\big)\Big),
\end{aligned}
\end{equation*}
where $I(\cdot)$ is taken from \eqref{e:a-index}, and where $\lambda^{j}_{i}=0$ if $i\in I(\ox(t)-\ou(t))\backslash I(x^j(t)-u^j(t))$ due to the clear inclusion $I(x^j(t)-u^j(t))\subset I(\ox(t)-\ou(t))$ for all $j=k,\ldots,\nu(k)$ and all large $k\in\N$.\vspace*{-0.05in}

Let us now show that the numerical sequence $\disp\left\{\sum^{\nu(k)}_{j=k}\al(k)_j\lm^j_i\right\}$ is bounded for all $i\in I(\ox(t)-\ou(t))$. Indeed, we deduce from \eqref{e:bound-grad} and \eqref{e:w-inverse-triangle-in} the validity of the estimates
$$
\begin{aligned}
\sum_{i\in I(\ox(t)-\ou(t))}&\sum^{\nu(k)}_{j=k}\al(k)_j\lm^j_i\le\dfrac{1}{M_1}\sum_{i\in I(\ox(t)-\ou(t))}\Big(\sum^{\nu(k)}_{j=k}\al(k)_j\lm^j_i\Big)\big\|\nabla g_i\big(\ox(t)-\ou(t)\big)\big\|\\
&\le\dfrac{\be}{M_1}\Big\|\sum_{i\in I(\ox(t)-\ou(t))}\Big(\sum^{\nu(k)}_{j=k}\al(k)_j\lm^j_i\Big)\nabla g_i\big(\ox(t)-\ou(t)\big)\Big\|\\
&\le\dfrac{\be}{M_1}\Big\|\sum_{i\in I(\ox(t)-\ou(t))}\sum^{\nu(k)}_{j=k}\al(k)_j\lm^j_i\nabla g_i\big(\ox(t)-\ou(t)\big)-\sum_{i\in I(\ox(t)-\ou(t))}\sum^{\nu(k)}_{j=k}\al(k)_j\lm^j_i\nabla g_i\big(x^j(t)-u^j(t)\big)\Big\|\\
&+\dfrac{\be}{M_1}\Big\|\sum_{i\in I(\ox(t)-\ou(t))}\sum^{\nu(k)}_{j=k}\al(k)_j\lm^j_i\nabla g_i\big(x^j(t)-u^j(t)\big)\Big\|\\
&\le\dfrac{\be}{M_1}\sum_{i\in I(\ox(t)-\ou(t))}\sum^{\nu(k)}_{j=k}\al(k)_j\lm^j_i\Big\|\nabla g_i\big(\ox(t)-\ou(t)\big)-\nabla g_i\big(x^j(t)-u^j(t)\big)\Big\|+\dfrac{\be\Tilde M}{M_1}\\
&\le\frac{1}{2}\sum_{i\in I(\ox(t)-\ou(t))}\sum^{\nu(k)}_{j=k}\al(k)_j\lm^j_i+\dfrac{\be\Tilde M}{M_1}
\end{aligned}
$$
for all $k$ sufficiently large, where $\Tilde M$ is an upper bound of $\disp\left\{\sum_{i\in I(\ox(t)-\ou(t))}\sum_{j=k}^{\nu(k)}\al(k)_j\lm^j_i\nabla g_i(x^j(t)-u^j(t))\right\}$. This justifies the boundedness of the sequence $\disp\left\{\sum^{\nu(k)}_{j=k}\al(k)_j\lm^j_i\right\}$. Thus there exists $\Tilde\al>0$ such that $\disp\left\|\sum^{\nu(k)}_{j=k}\al(k)_j\lm^j_i\right\|\le\Tilde\al$ and $\disp\sum^{\nu(k)}_{j=k}\al(k)_j\lm^j_i \to \be_i$ as $k\to\infty$ along some subsequence, with no relabeling.

Next we verify that $\sum_{j=k}^{\nu(k)}\al(k)_j\lm^j_i\nabla g_i(x^j(t)-u^j(t))\to\be_i\nabla g_i(\ox(t)-\ou(t))$ as $k\to\infty$. Observe that
$$
\begin{aligned}
\Big\|\sum_{j=k}^{\nu(k)}&\al(k)_j\lm^j_i\nabla g_i\big(x^j(t)-u^j(t)\big)-\be_i\nabla g_i\big(\ox(t)-\ou(t)\big)\Big\|\\
&\le\Big\|\sum^{\nu(k)}_{j=k}\al(k)_j\lm^j_i\nabla g_i\big(x^j(t)-u^j(t)\big)-\sum^{\nu(k)}_{j=k}\al(k)_j\lm^j_i\nabla g_i\big(\ox(t)-\ou(t)\big)\Big\|\\
&+\Big\|\sum^{\nu(k)}_{j=k}\al(k)_j\lm^j_i\nabla g_i\big(\ox(t)-\ou(t)\big)-\be_i\nabla g_i\big(\ox(t)-\ou(t)\big)\Big\|\\
&\le\sum^{\nu(k)}_{j=k}\al(k)_j\lm^j_i\Big\|\nabla g_i\big(x^j(t)-u^j(t)\big)-\nabla g_i\big(\ox(t)-\ou(t)\big)\Big\|+M_2\Big\|\sum^{\nu(k)}_{j=k}\al(k)_j\lm^j_i-\beta_i\Big\|\\
&\le\Tilde\al\ve+M_2\ve\;\mbox{ for all large }\;k\in\N,
\end{aligned}
$$
where $M_2$ is taken in (H2) while $\ve>0$ is an arbitrary small number with
$$
\max\Big\{\big\|\nabla g_i\big(x^j(t)-u^j(t)\big)-\nabla g_i\big(\ox(t)-\ou(t)\big)\big\|,\Big\|\sum^{\nu(k)}_{j=k}\al(k)_j\lm^j_i-\beta_i\Big\|\Big\}\le\ve
$$
whenever $k\le j\le\nu(k)$ and $k$ is large enough. Thus we have
\begin{equation*}
\begin{aligned}
-\dot{\ox}(t)-&f\big(\ox(t),\ta(t)\big)=\lim_{k\to\infty}\Big(-\sum_{i\in I(\ox(t)-\ou(t))}\sum_{j=k}^{\nu(k)}\al(k)_j\lm^j_i\nabla g_i\big(x^j(t)-u^j(t)\big)\Big)\\
&=-\sum_{i\in I(\ox(t)-\ou(t))}\be_i\nabla g_i\big(\ox(t)-\ou(t)\big)\in N\big(\ox(t)-\ou(t);C\big)\;\mbox{ a.e. }\;t\in[0,T],
\end{aligned}
\end{equation*}
which verifies by \eqref{e:F-rep} that $\oz(\cdot)$ satisfies the differential inclusion \eqref{e:diff-incl} and hence the constraints in \eqref{e:mc}.\vspace*{-0.05in}

To justify the optimality of the triple $\oz(\cdot)$ in $(P)$, it remains to show that
\begin{eqnarray*}
J[\ox,\ou,\oa]\le\liminf_{k\to\infty}J[x^k,u^k,a^k]
\end{eqnarray*}
for the Bolza functional \eqref{e:Bolza}. But this is a clear consequence of the aforementioned Mazur theorem due the imposed convexity of the integrand $\ell$ in the velocity variables and the Lebesgue dominated convergence theorem for passing to the limit under the integral sign. $\h$\vspace*{-0.02in}

Note that the convexity of the running cost $\ell$ in \eqref{e:Bolza} is not among the standing assumptions of the paper and is not needed for deriving our main results on necessary optimality conditions in $(P)$. Such conditions established below address the so-called ``relaxed intermediate local minimizers" introduced in \cite{m95} and then studied in many publications. To recall this notion, we first consider the relaxed optimal control problem for $(P)$ following the Bogoluybov-Young relaxation/convexification procedure, which has been well understood in the classical calculus of variations and optimal control; see, e.g., \cite{m-book2} with the references therein and \cite{cmr,dfm,et,tols} for more recent results in this direction. To this end, denote by $\hat\ell_F(t,x,u,a,\dot{x},\dot{u},\dot{a})$ the convexification (the largest l.s.c.\ convex function majorized by $\ell(t,x,u,a,\cdot,\cdot,\cdot)$) of the running cost in \eqref{e:Bolza} on the set $F(x,u,a)$ from \eqref{e:vel-map} with respect to the velocity variables $(\dot{x},\dot{u},\dot{a})$ and put $\hat\ell:=\infty$ at points out of $F(x,u,a)$. Define the {\em relaxed optimal control problem} $(R)$ by
\begin{equation}\label{e:relaxed-problem}
\mbox{minimize}\;\;\hat J[z]=\hat J[x,u,a]:=\ph\big(x(T)\big)+\int^T_0\hat\ell_F\big(t,x(t),u(t),a(t),\dot{x}(t),\dot{u}(t),\dot{a}(t)\big)dt
\end{equation}
over all $z(\cdot)=(x(\cdot),u(\cdot),a(\cdot))\in W^{1,2}([0,T])$ satisfying \eqref{e:mset}, \eqref{e:ec}, and \eqref{e:diff-incl}. Besides the obvious case of integrands that are convex in velocity variables, there are broad classes of variational and control problems over continuous-time intervals where optimal values of the cost functionals in the (nonconvex) original and relaxed problems agree; it is known as ``relaxation stability." This is due to some``hidden convexity" for such problems (we refer again to \cite[Chapter~6]{m-book2} and the commentaries therein), which allows us, in particular, to verify the relaxation stability of nonconvex Bolza problems for Lipschitzian differential inclusions and also for those satisfying a certain one-sided Lipschitzian condition \cite{dfm}. Unfortunately, neither of the aforementioned conditions is fulfilled for the controlled sweeping process under consideration. For a sweeping process over prox-regular moving sets with controls only in additive perturbations, the relaxation stability follows from the result by Edmond and Thibault \cite[Theorem~2]{et}; see also \cite{cmr} for similar relaxation results concerning BV solutions. It seems that the closest to our setting is a recent result by Tolstonogov \cite[Theorem~4.2]{tols}, which establishes the relaxation stability for a sweeping process over variable convex moving sets involved in optimization together with controls in perturbations.\vspace*{-0.02in}

Following \cite{m95,m-book2}, we say that $\oz(\cdot)\in W^{1,2}([0,T];\R^{2n+d})$ is a {\em relaxed intermediate local minimizer} (r.i.l.m.) of rank 2 for $(P)$ if it is feasible to $(P)$, $J[\oz]=\hat J[\oz]$, and there is $\ve>0$ such that $J[\oz]\le J[z]$ for any feasible solution $z(\cdot)$ to $(P)$ satisfying the conditions
\begin{equation}\label{e:rilm}
\|z(t)-\oz(t)\|<\ve\;\mbox{ whenever }\;t\in[0,T]\;\mbox{ and }\;\int^T_0\|\dot{z}(t)-\dot{\oz}(t)\|^2dt<\ve.
\end{equation}
Since we do not consider in this paper relaxed intermediate local minimizers of any other rank but 2, we skip mentioning the rank in what follows. The name ``intermediate" comes from the fact that the introduced notion clearly lies (strictly) between the conventional notions of weak and strong local minimizers in the calculus of variations and optimal control; see \cite{m-book2} for more discussions.\vspace*{-0.05in}

Note also that in the case of $J[\oz]=\hat J[\oz]$ (in particular, if $(P)$ has the property of relaxation stability), there is no difference between relaxed intermediate local minimizers and merely intermediate local minimizers (without relaxation), which were also defined in \cite{m95}. Thus we can treat r.i.l.m. as a {\em local} version of relaxation stability. It can be distilled from the proofs of \cite[Theorem~2]{et} and \cite[Theorem~4.2]{tols} that no relaxation is needed provided that $\oz(\cdot)$ ia a {\em strong} local minimizer of $(P)$ and, in addition to our standing assumptions, either controls are presented only in perturbations, or the set $C$ in \eqref{e:mset} is convex.\vspace*{-0.2in}

\section{Discrete Approximation of Intermediate Local Minimizers}
\setcounter{equation}{0}\vspace*{-0.1in}

In this section we continue with developing the method of discrete approximation, while now paying our main attention not to constructing such approximations of {\em any feasible} solution to the constrained system \eqref{e:diff-incl}, \eqref{e:mset}, and \eqref{e:ec} as in Section~3 but to the given {\em local optimal solution} $\oz(\cdot)=(\ox(\cdot),\ou(\cdot),\oa(\cdot))$ (in the {\em r.i.l.m.} sense) of this system with respect to the cost functional \eqref{e:Bolza}, which we can write as
\begin{equation}\label{e:Bolzaz}
\mbox{minimize}\;\;J[z]=\ph\big(z(T)\big)+\int^T_0\ell\big(t,z(t),\dot{z}(t)\big)dt
\end{equation}
with $\ph(z):=\ph(x)$. It means that our goal here is to construct discrete approximations of the entire problem $(P)$ including its cost functional and the given r.i.l.m.\ for which we aim subsequently to derive necessary optimality conditions by employing the discrete approximation method.\vspace*{-0.02in}

For any fixed $k\in\N$ we define the {\em discrete sweeping control problem} $(P_k)$ as follows: minimize
\begin{equation}
\label{e:discrete-cost}
\begin{aligned}
J_k[z^k]:&=\ph(x^k_k)+h_k\sum^{k-1}_{j=0}\ell\left(t^k_j,x^k_j,u^k_j,a^k_j,\dfrac{x^k_{j+1}-x^k_j}{h_k},\dfrac{u^k_{j+1}-u^k_j}{h_k},\dfrac{a^k_{j+1}-a^k_j}{h_k}\right)\\
&+\sum^{k-1}_{j=0}\int^{t^k_{j+1}}_{t^k_j}\left(\bigg\|\dfrac{x^k_{j+1}-x^k_j}{h_k}-\dot{\ox}(t)\bigg\|^2+\bigg\|\dfrac{u^k_{j+1}-u^k_j}{h_k}-\dot{\ou}(t)\bigg\|^2+\bigg\|\dfrac{a^k_{j+1}-a^k_j}{h_k}-\dot{\oa}(t)\bigg\|^2\right)\\
&+\dist^2\left(\bigg\|\dfrac{u^k_1-u^k_0}{h_k}\bigg\|;(-\infty,\Tilde\mu]\right)+\dist^2\left(\sum^{k-2}_{j=0}\bigg\|\dfrac{u^k_{j+2}-2u^k_{j+1}+u^k_j}{h_k}\bigg\|;(-\infty,\Tilde\mu]\right)
\end{aligned}
\end{equation}
over elements $z^k:=(x^k_0,x^k_1,\ldots,x^k_k,u^k_0,u^k_1,\ldots,u^k_{k-1},a^k_0,a^k_1,\ldots,a^k_{k-1})$ satisfying the constraints
\begin{equation}\label{e:dc1}
x^k_{j+1}\in x^k_j-h_kF(x^k_j,u^k_j,a^k_j)\;\;\mbox{for}\;\;j=0,\ldots,k-1\;\;\mbox{with}\;\;(x^k_0,u^k_0,a^k_0)=\big(x_0,\ou(0),\oa(0)\big),
\end{equation}
\begin{equation}\label{e:dc2}
g_i(x^k_k-u^k_k)\ge 0\;\;\mbox{for}\;\;i=1,\ldots,m,
\end{equation}
\begin{equation}\label{e:dc3}
r_1-\ve_k\le\|u^k_j\|\le r_2+\ve_k\;\;\mbox{for}\;\;j=0,\ldots,k,
\end{equation}
\begin{equation}\label{e:dc5}
\|(x^k_j,u^k_j,a^k_j)-(\ox(t^k_j),\ou(t^k_j),\oa(t^k_j))\|\le\ve/2\;\;\mbox{for}\;\;j=0,\ldots,k-1,
\end{equation}
\begin{equation}\label{e:dc6}
\sum^{k-1}_{j=0}\int^{t^k_{j+1}}_{t^k_j}\left(\bigg\|\dfrac{x^k_{j+1}-x^k_j}{h_k}-\dot{\ox}(t)\bigg\|^2+\bigg\|\dfrac{u^k_{j+1}-u^k_j}{h_k}-\dot{\ou}(t)\bigg\|^2+
\bigg\|\dfrac{a^k_{j+1}-a^k_j}{h_k}-\dot{\oa}(t)\bigg\|^2\right)dt\le\dfrac{\ve}{2},
\end{equation}
\begin{equation}\label{e:dc7}
\bigg\|\dfrac{u^k_1-u^k_0}{t^k_1-t^k_0}\bigg\|\le\Tilde\mu+1\;\;\mbox{and}\;\;\sum^{k-2}_{j=0}\bigg\|\dfrac{u^k_{j+2}-2u^k_{j+1}+u^k_j}{h_k}\bigg\|\le\Tilde\mu+1,
\end{equation}
where $\ve$ is taken from \eqref{e:rilm} while $\ve_k$ and $\Tilde\mu$ are taken from Theorem~\ref{Th:DA} applied to the given r.i.l.m.\ $\oz(\cdot)$.\vspace*{-0.02in}

To study $\oz(\cdot)$ via the method of discrete approximations, we have to verify first that all the problems $(P_k)$ for each $k\in\N$ sufficiently large admit optimal solutions.\vspace*{-0.1in}

\begin{proposition}{\bf (existence of optimal solutions to discrete sweeping control problems ).}\label{Pr:existence-DA} Under the standing assumptions {\rm (H1)--(H4)} holding along the r.i.l.m.\ $\oz(\cdot)$ each problem $(P_k)$  for large $k\in\N$ admits an optimal solution $\oz^k(\cdot)$.
\end{proposition}\vspace*{-0.1in}
{\bf Proof.} It follows directly from Theorem~\ref{Th:DA} and the construction  of $(P_k)$ along the given r.i.l.m.\ $\oz(\cdot)$ that the set of feasible solutions to $(P_k)$ is nonempty for all large $k\in\N$. Furthermore, the imposed constraints \eqref{e:dc3}-\eqref{e:dc6} ensure that this set is bounded for each $k$. Thus the existence of optimal solutions to finite-dimensional problems $(P_k)$ is ensures by the classical the Weierstrass existence theorem provided that the feasible solution set to each problem $(P_k)$ is closed. To check it, take a sequence $z^\nu(\cdot)=z^\nu:=(x^\nu_0,\ldots,x^\nu_k,u^\nu_0,\ldots,u^\nu_{k-1},a^\nu_0,\ldots,a^\nu_{k-1})$ of feasible solutions for $(P_k)$ converging to some $z(\cdot)=z:=(x_0,\ldots,x_k,u_0,\ldots,u_{k-1},a_0,\ldots,a_{k-1})$ as $\nu\to\infty$. We need to show that $z$ is feasible to $(P_k)$. Observing that
$g_i(x_j-u_j)=\disp\lim_{\nu\to\infty}g_i(x^\nu_j-u^\nu_j)\ge0$ for all $i=1,\ldots,m$ and $j=0,\ldots,k-1$ gives us that $x_j-u_j\in C$ for all $j=0,\ldots,k-1$. It is not hard to see that $I(x^\nu_j-u^\nu_j)\subset I(x_j-u_j)$ for $\nu\in\N$ sufficiently large. Taking \eqref{e:vel-map} and \eqref{e:diff-incl} into account, gives us for all such indices $j$ that
$$
\dfrac{x^\nu_{j+1}-x^\nu_j}{-h_k}-f(x^\nu_j,a^\nu_j)\in N_C(x^\nu_j-u^\nu_j),
$$
which implies therefore by passing to the limit as as $\nu\to\infty$ that
$$
\dfrac{x^\nu_{j+1}-x^\nu_j}{-h_k}-f(x^\nu_j,a^\nu_j)\to\dfrac{x_{j+1}-x_j}{-h_k}-f(x_j,a_j)\;\;\mbox{and}\;\;x^\nu_j-u^\nu_j\to x_j-u_j.
$$
This allows us to arrive at at the inclusions
$$
\dfrac{x_{j+1}-x_j}{-h_k}-f(x_j,a_j)\in\Limsup_{x-u\to x_j-u_j}N_C(x-u)=N_C(x_j-u_j),
$$
ensuring that $x_{j+1}-x_j\in F(x_j,u_j,a_j)$ for all $j=0,\ldots,k-1$ and thus completing the proof.$\h$\vspace*{-0.02in}

The next key result makes a bridge between optimal solutions to the control problems $(P)$ and $(P_k)$ by showing on one hand that the optimal solutions $\oz^k(\cdot)$ to $(P_k)$ are {\em approximately optimal/suboptimal} solutions to $(P)$ and that necessary optimality conditions for $\oz^k$ can be treated as ``almost optimality conditions" for $\oz(\cdot)$. On the other hand, the necessary optimality conditions for the discrete solutions $\oz^k(\cdot)$ obtained below will serve as the basis to derive the exact necessary optimality conditions for $\oz(\cdot)$ by passing to the limit in their relationships as $k\to\infty$.\vspace*{-0.1in}

\begin{theorem}{\bf (strong $W^{1,2}$ convergence of discrete optimal solutions).}\label{Th:strong-DA-ilm} Suppose that all the standing assumptions {\rm(H1)--(H4)} and those of Theorem~{\rm\ref{Th:DA}} hold along the given r.i.l.m.\ $\oz(\cdot)$ of problem $(P)$. Assume in addition that the terminal cost $\ph$ is continuous at $\ox(T)$, that the running cost $\ell$  in \eqref{e:Bolzaz}) is continuous at $(t,\oz(t),\dot{\oz}(t))$ for a.e.\ $t\in[0,T]$, and that $\ell(\cdot,z,\dot{z})$ is uniformly majorized  around $\oz(\cdot)$ by a summable function on $[0,T]$. Then any sequence of optimal solutions $\oz^k(\cdot)=(\ox^k(\cdot),\ou^k(\cdot),\oa^k(\cdot))$ of $(P_k)$, piecewise linearly extended to the whole interval $[0,T]$, converges to $\oz(\cdot)$ in the strong topology of $W^{1,2}([0,T];\R^{2n+d})$ with the validity of the estimates
\begin{equation}\label{e:discrete-estimate}
\bigg\|\dfrac{\ou^k_1-\ou^k_0}{h_k}\bigg\|\le\Tilde\mu\;\mbox{ and }\;\limsup_{k\to\infty}\sum^{k-2}_{j=0}\bigg\|\dfrac{\ou^k_{j+2}-2\ou^k_{j+1}+\ou^k_j}{h_k}\bigg\|\le\Tilde\mu,
\end{equation}
where the number $\Tilde\mu>0$ is taken from Theorem~{\rm\ref{Th:DA}}.
\end{theorem}\vspace*{-0.1in}
{\bf Proof.} Take any sequence $\{\oz^k(\cdot)\}$ of optimal solutions to $(P_k)$, the existence of which is ensured by Proposition~\ref{Pr:existence-DA}, and then each $\oz^k(\cdot)$ piecewise linearly to the continuous-time interval $[0,T]$. All the statements of the theorem follow from the following limiting equality:
\begin{equation}\label{e:52}
\begin{aligned}
\lim_{k\to\infty}\int^{T}_{0}&\bigg(\left\|\dot{\ox}(t)-\dot{{\ox}}^{k}(t)\right\|^{2}+\left\|\dot{\ou}(t)-\dot{{\ou}}^{k}(t)\right\|^{2}+
\left\|\dot{\oa}(t)-\dot{{\oa}}^{k}(t)\right\|^{2}\bigg)dt\\
+&\Big\|\frac{\ox^k_1-\ox^k_0}{h_k}-\dot\ox(0)\Big\|^2+\dist^2\bigg(\bigg\|\dfrac{\ou^k_1-\ou^k_0}{h_k}\bigg\|;\big(-\infty,\Tilde{\mu}\big]\bigg)\\
+&\dist^2\bigg(\sum^{k-2}_{j=0}\bigg\|\dfrac{\ou^k_{j+2}-2\ou^k_{j+1}+\ou^k_j}{h_k}\bigg\|;\big(-\infty,\Tilde{\mu}\big]\bigg)=0.
\end{aligned}
\end{equation}
Arguing by contradiction, suppose that \eqref{e:52} fails and thus find a subsequence of $k\in\N$ along which the limit in \eqref{e:52} equals to some $c>0$. The weak compactness of the unit ball in $L^{2}([0,T];\R^{2n+d})$ yields the existence of a triple $(v(\cdot),w(\cdot),q(\cdot))\in L^{2}([0,T];\R^{2n+d})$ such that
$$
\big(\dot{{\ox}}^{k}(\cdot),\dot{{\ou}}^{k}(\cdot),\dot{{\oa}}^{k}(\cdot)\big)\to\big(v(\cdot),w(\cdot),q(\cdot)\big)\;\mbox{ weakly in }\;L^{2}([0,T];\R^{2n+d})
$$
along a subsequence of optimal velocities $\{\dot{\oz}^k(\cdot)\}$. It is clear that $\dot{\tz}(t)=(v(t),w(t),q(t))$ a.e. on $[0,T]$ for the absolutely continuous triple defined by
$$
\tz(t):=\big(x_{0},\bar{u}(0),\bar{a}(0)\big)+\int^{t}_{0}(v(s),w(s),q(s)\big)ds\;\mbox{ for all }\;t\in[0,T].
$$
The latter implies that $\dot\oz^k(\cdot)\to\dot{\tz}(\cdot)=(\dot{\tx}(\cdot),\dot{\tu}(\cdot),\dot{\ta}(\cdot))$ weakly in $L^{2}([0,T];\R^{2n+d})$, and hence it yields  $\tz(\cdot)\in W^{1,2}([0,T];\R^{2n+d})$. Similarly to the proof of Theorem~\ref{Th:existence} we verify that $\tz(\cdot)$ satisfies inclusion \eqref{e:diff-incl} with the constraints in \eqref{e:ec} and \eqref{e:mc}. The rest of the proof of this theorem follows the lines in the proof of \cite[Theorem~5.1]{cm1}, which show that $\tz(\cdot)$ is feasible to the relaxed problem $(R)$, belongs to the selected $W^{1,2}$ neighborhood of $\oz(\cdot)$, and gives a smaller value to \eqref{e:relaxed-problem} than $\oz(\cdot)$ with $\hat J[\oz]=J[\oz]$. Thus the assumed failure of \eqref{e:52} contradicts the choice of $\oz(\cdot)$ as a r.i.l.m.\ in $(P)$, and we are done. $\h$\vspace*{-0.2in}

\section{Second-Order Subdifferential Computations}
\setcounter{equation}{0}\vspace*{-0.1in}

Optimization problems $(P)$ and $(P_k)$ are intrinsically {\em nonsmooth} and {\em nonconvex}, even for smooth and/or convex terminal and running costs. The unavoidable source of nonsmoothness and nonconvexity comes from  the sweeping differential inclusion \eqref{e:diff-incl} and its discrete approximations \eqref{e:dc1}, which constitute nonconvex geometric constraints of the {\em graphical} type. Furthermore, the first-order normal cone (subdifferential, variational structure) of the sweeping inclusions \eqref{e:diff-incl} and \eqref{e:dc1} calls for appropriate {\em second-order} subdifferential constructions to derive and analyze optimality conditions for their solutions. In this section we recall the corresponding generalized differential constructions and present the results of their computations in terms of the initial problem data that play a significant role in what follows.\vspace*{-0.02in}

Given a set-valued mapping $F\colon\R^n\tto\R^m$, we always assume that its graph
$$
\gph F:=\big\{(x,y)\in\R^n\times\R^m\big|\;y\in F(x)\big\}
$$
is locally closed around the reference point $(\ox,\oy)\in\gph F$ and define its {\em coderivative} of $F$ at this point via the (limiting) normal cone \eqref{e:Mor-nc} by
\begin{equation}\label{e:cor}
D^*F(\ox,\oy)(u):=\big\{v\in\R^n\;|\;(v,-u)\in N\big((\ox,\oy);\gph F\big)\big\}\;\mbox{for all }\;u\in\R^m.
\end{equation}
When $F\colon\R^n\to\R^m$ is single-valued (then $\oy=F(\ox)$ is omitted in the coderivative notation) and continuously differentiable $({\cal C}^1$-smooth) around $\ox$, we have the representation
$$
D^*F(\ox)(u)=\big\{\nabla F(\ox)^*u\big\},\quad u\in\R^m,
$$
via the adjoint/transposed Jacobian matrix $\nabla F(\ox)^*$. The corresponding (first-order) {\em subdifferential} of an l.s.c.\  function $\ph\colon\R^n\to\oR$ at $\ox\in\dom\ph:=\{x\in\R^n|\;\phi(x)<\infty\}$ can be defined geometrically
\begin{equation}\label{e:sub}
\partial\phi(\ox):=\big\{v\in\R^n\big|\;(v,-1)\in N\big((\ox,\ph(\ox)\big);\epi\phi)\big\}
\end{equation}
via the normal cone \eqref{e:Mor-nc} of its epigraph $\epi\phi:=\{(x,\al)\in\R^{n+1}|\;\al\ge\phi(x)\}$ while admitting, together with the coderivative \eqref{e:cor} equivalent analytical representations and---the crucial issue--satisfy  comprehensive {\em calculus rules} despite the nonconvexity of their values; see \cite{m-book1,rw} with the references therein.\vspace*{-0.03in}

Now we turn to the {\em second-order subdifferential}/generalized Hessian of $\phi\colon\R^n\to\oR$ at $\ox\in\dom\phi$ relative to $\ov\in\partial\phi(\ox)$, which plays an underlying role in this paper and is defined, following the ``dual derivative-of-derivative" scheme \cite{m92}, as the coderivative of the first-order subdifferential of $\phi$ by
\begin{eqnarray}\label{2nd}
\partial^2\phi(\ox,\ov)(u):=\big(D^*\partial\phi\big)(\ox,\ov)(u),\quad u\in\R^n.
\end{eqnarray}
If $\phi$ is ${\cal C}^2$-smooth around $\ox$, then \eqref{2nd} reduces to the (symmetric) Hessian matrix $\partial^2\phi(\ox)(u)=\{\nabla^2\phi(\ox)u\}$ for all $u\in\R^n$, while in general it is a positively homogeneous set-valued mapping of $u$ satisfying well-developed {\em second-order subdifferential calculus}; see \cite{m-book2,mr}.\vspace*{-0.02in}

To present next the second-order computations needed in what follows, we recall the definition of calmness, which is a weak ``one-point" stability property that has been well-understood by now in variational analysis and optimization; see \cite{rw}, \cite{hos} and the references therein. A set-valued mapping $M\colon\R^s\tto\R^q$ is {\em calm} at $(\bar\vt,\bar q)\in\gph M$ if there are positive numbers $\mu$ and $\eta$ such that
\begin{equation}\label{e:calm}
M(\vt)\cap(\bar q+\eta\B)\subset M(\bar\vt)+\mu\|\vt-\bar\vt\|\B\;\mbox{ whenever }\;\vt\in\bar\vt+\eta\B.
\end{equation}
Due to the normal cone description of the sweeping process, in this paper we employ the second-order subdifferential \eqref{2nd} just for the set indicator function $\dd_\O(x)$ of $\O\subset\R^n$ that is equal to $0$ if $x\in\O$ and to $\infty$ if $x\notin\O$. The second-order subdifferential of the indicator function clearly reduces to $D^*N_\O$, where we use the notation $N_\O(x):=N(x;\O)$ for convenience. In this case the following upper estimates and exact computations of $D^*N_\O$ can be deduced from \cite{hos} and the previous developments mentioned therein.\vspace*{-0.1in}

\begin{proposition}{\bf (coderivative of the normal cone mapping to inequality constraints).}\label{Th:co-nc} Consider the set $\O:=\{x\in\R^n|\;g_i(x)\ge 0\}$ defined by the ${\cal C}^2$-smooth functions $g=(g_1,\ldots,g_m)\colon\R^n\to\R^m$ around $\ox\in\O$ so that the vectors $\nabla g_1(\ox),\ldots\nabla g_m(\ox)$ are positively linearly independent, which amounts to saying that the Mangasarian-Fromovitz constraint qualification $($MFCQ$)$ is satisfied at $\ox$. Given a normal $\ov\in N_\O(\ox)$, suppose in addition that the multifunction $M\colon\R^{2m}\tto\R^{n+m}$ defined by
\begin{eqnarray}\label{calm}
M(\vTh):=\big\{(x,\lm)\big|\;\big(-g(x),\lm\big)+\vTh\in\gph N_{\R^m_-}\big\}
\end{eqnarray}
is calm \eqref{e:calm} at $(0,\ox,\bar\lm)$ for all $\bar\lm=(\bar\lm_1,\ldots,\bar\lm_m)\ge 0$ satisfying the equation $-\nabla g(\ox)^*\bar\lm=\bar v$. Then we have the second-order upper estimate
\begin{eqnarray*}
D^*N_\O(\ox,\bar v)(u)\subset\bigcup_{\bar\lm\ge 0,-\nabla g(\ox)\bar\lm=\bar v}\left\{\left(-\sum^m_{i=1}\bar\lm_i\nabla^2g_i(\ox)\right)u-\nabla g(\ox)^*D^*N_{\R^m_-}\big(-g(\ox),\bar\lm\big)\big(-\nabla g(\ox)u\big)\right\}.
\end{eqnarray*}
Strengthening the calmness assumption by the full rank of the Jacobian $\nabla g(\ox)$ gives us the precise formula
\begin{equation*}
D^*N_\O(\ox,\bar v)(u)=\left(-\sum^m_{i=1}\bar\lm_i\nabla^2g_i(\ox)\right)u-\nabla g(\ox)^*D^*N_{\R^m_-}\big(-g(\ox),\bar\lm\big)\big(-\nabla g(\ox)u\big),
\end{equation*}
where $\bar\lm\ge 0$ is a unique solution to the equation $-\nabla g(\ox)^*\bar\lm=\bar v$. Furthermore, the coderivative of the normal cone mapping generated by the nonnegative orthant $\R^m_-$ above is computed by
\begin{eqnarray}\label{e:cor-nc-o}
D^*N_{\R^m_-}(x,v)(y)=\left\{\begin{array}{ll}
\emp&\mbox{if }\,\exists\;i\;\mbox{ with }\;v_iy_i\not=0,\\
\{\gg|\;\gg_i=0\;\forall\;i\in I_1(y),\;\gg_i\ge 0\;\forall\;i\in I_2(y)\big\}&\mbox{otherwise}
\end{array}\right.
\end{eqnarray}
whenever $(x,v)\in\gph N_{\R^m_-}$ with the index subsets in \eqref{e:cor-nc-o} defined by
\begin{equation}\label{e:indexes}
I_1(y):=\big\{i\big|\;x_i<0\big\}\cup\big\{i\big|\;v_i=0,\;y_i<0\big\},\quad I_2(y):=\big\{i\big|\;x_i=0,\;v_i=0,\;y_i>0\big\}.
\end{equation}
\end{proposition}\vspace*{-0.1in}
{\bf Proof.} Compare \cite[Theorem~3.3]{hos}, its proof, and the further references therein. $\h$\vspace*{-0.02in}

The next result present the crucial second-order computations of the coderivative of the sweeping process under consideration entirely in terms of its given data.\vspace*{-0.1in}

\begin{theorem}{\bf (second-order computations for the sweeping process).}\label{Th:co-cal} Consider the set-valued mapping $F$ associated with the sweeping process \eqref{e:5} by \eqref{e:vel-map}, where the nonconvex set $C$ is taken from \eqref{e:mset}, and where the perturbation mapping that $f$ is ${\cal C}^1$-smooth. Given $x,u\in\R^n$ with $x-u\in C$ as well as $w\in N_C(x-u)$ and $a\in\R^d$, suppose that the vectors $\nabla g_1(x-u),\ldots,\nabla g_m(x-u)$ are positively linearly independent and that the multifunction $M$ from \eqref{calm} is calm \eqref{e:calm} at $(0,x-u,\lm)$ for all $\lm=(\lm_1,\ldots,\lm_m)\ge 0$ satisfying the equation $-\nabla g(x-u)^*\lm=w-f(x,a)$. Then we have the upper estimate
\begin{eqnarray}\label{e:co-upper-est1}
\begin{array}{ll}
D^*F(x,u,a,w)(y)\subset\disp\bigcup_{\lm\ge 0,-\nabla g(x-u)\lm= w-f(x,a)}\bigg\{\bigg(\nabla_xf(x,a)^*y-\bigg(\sum^m_{i=1}\lm\nabla^2g_i(x-u)\bigg)y-\nabla g(x-u)^*\gg,\\ \disp\bigg(\sum^m_{i=1}\lm\nabla^2g_i(x-u)\bigg)y+\nabla g(x-u)^*\gg,\nabla_af(x,a)^*y\bigg)\bigg\}\;\mbox{ for all }\;y\in\dom D^*N_C\big(x-u,w-f(x,a)\big),
\end{array}
\end{eqnarray}
where the coderivative domain is satisfied the inclusion
\begin{equation}\label{e:dom-co-est}
\begin{aligned}
\dom D^*N_C(x-u,w-f(x,a))\subset&\big\{y|\;\exists\,\lm\ge0\;\mbox{ such that }\;-\nabla g(x-u)\lm=w-f(x,a),\\
&\lm_i\la\nabla g_i(x-u),y\ra=0\;\mbox{ for }\;i=1,\ldots,m\big\},
\end{aligned}
\end{equation}
and where we have in \eqref{e:co-upper-est1} that $\gg_i=0$ if either $g_i(x-u)>0$ or $\lm_i=0$ and $\la\nabla g_i(x-u),y\ra>0$, and that $\gg_i\ge 0$ if $g_i(x-u)=0,\lm_i=0$, and $\la\nabla g_i(x-u),y\ra<0$.\\[1ex]
Furthermore, replacing the calmness of \eqref{calm} by the stronger assumptions on the full rank of the Jacobian matrix $\nabla g(x-u)$ $($which is actually the classical LICQ--linear independence constraint qualification$)$ ensures that the equalities hold in both inclusions \eqref{e:co-upper-est1} and \eqref{e:co-upper-est1} with the collection of nonnegative multipliers $\lm=(\lm_1,\ldots,\lm_m)\ge 0$ uniquely determined by the equation $-\nabla g(x-u)^*\lm=w-f(x,a)$.
\end{theorem}\vspace*{-0.1in}
{\bf Proof.} Consider the mappings $G(x,u,a):=N_C(x-u)$ and $\Tilde f(x,u,a):=f(x,a)$. It follows from the coderivative sum rule in \cite[Theorem~1.62]{m-book1} that
$$
z^*\in\nabla\Tilde f(x,u,a)^*y+D^*G\big(x,u,a,w-f(x,a)\big)(y)
$$
for any $y\in\dom D^*N_C(x-u,w-f(x,a))$ and $z^*\in D^*F(x,u,a,w)(y)$. Observe further that
$$
G(x,u,a)=N_C\circ\Tilde g(x,u,a)\;\mbox{ with }\;\Tilde g(x,u,a):=x-u,
$$
where the Jacobian of latter mapping is obviously of full rank. It follows from the coderivative chain rule of \cite[Theorem~1.66]{m-book1} applied to the above composition that
\begin{equation}\label{e:cf}
z^*\in\nabla\Tilde f(x,u,a)^*y+\nabla\Tilde g(x,u,a)^*D^*N_C\big(x-u,w-f(x,a)\big)(y).
\end{equation}
Substituting now into \eqref{e:cf} the corresponding results of Proposition~\ref{Th:co-nc} and taking into account the structure of the mapping $\Tilde f$ in \eqref{e:cf} give us all the statements claimed in the theorem. $\h$\vspace*{-0.02in}

Note that for the linear case of $g_i(x)=-\la x^*_i,x\ra$, corresponding to the polyhedral sweeping process, the upper estimate in \eqref{e:co-upper-est1} reduces to our previous computations in \cite[Theorem~6.1]{cm1}.\vspace*{-0.2in}

\section{Necessary Optimality Conditions for Discrete-Time Problems}
\setcounter{equation}{0}\vspace*{-0.1in}

This section is devoted to deriving necessary optimality conditions for local optimal solution to the discrete-time control problems $(P_k)$, for each fixed $k\in\N$. First we establish, under minimal assumptions, necessary optimality conditions in the extended Euler-Lagrange form for a general class of problems $(P_k)$ with an arbitrary discrete velocity map $F$ via its coderivative by reducing such problems to nonsmooth mathematical programming with many geometric constraints of the graphical type. Then we exploit the special normal cone structure \eqref{e:vel-map} of $F$ to obtain optimality conditions for the discrete sweeping control problems $(P_k)$ defined in Section~5 expressed entirely via the given data by using the second-order computations of Section~6. Due to the approximation results of Section~5, the optimality conditions for $(P_k)$ obtained in this way can be treated as {\em suboptimality} conditions for the given r.i.l.m.\ of the original sweeping control problem $(P)$, while our main goal in what follows is to derive necessary optimality conditions for such local minimizers of $(P)$ by passing to the limit from those in discrete approximations.\vspace*{-0.1in}

\begin{theorem}{\bf (necessary conditions of the Euler-Lagrange type for discrete-time optimal control).}\label{Th:EL-DA} Let $\oz^k=(x_0,\ox^k_1,\ldots,\ox^k_k,\ou^k_0,\ldots,\ou^k_k,\oa^k_0,\ldots,\oa^k_k)$ be a local optimal solution to problem $(P_k)$ for whenever $k\in\N$, where $F$ is an arbitrary closed-graph mapping, and where $\ph$ and $\ell_t:=\ell(t,\cdot,\cdot$ are locally Lipschitzian around the optimal points for any $t\in\Delta_k$.
Then there exist dual elements
$\al^k=(\al^k_1,\ldots,\al^k_m)\in\R^m_+,\;\xi^{1k}=(\xi^{1k}_0,\ldots,\xi^{1k}_k)\in\R^{k+1}_+,\;\xi^{2k}=(\xi^{2k}_0,\ldots,\xi^{2k}_k)\in\R^{k+1}_-$, and $p^k_j=(p^{xk}_j,p^{uk}_j,p^{ak}_j)\in\R^n\times\R^n\times\R^d$ for $j=0,\ldots,k$ satisfying the conditions
\begin{equation}\label{e:nontri-da}
\lm^k+\|\al^k\|+\|\be^k\|+\|\xi^{1k}+\xi^{2k}\|+\sum^{k-1}_{j=0}\|p^{xk}_j\|+\|p^{uk}_0\|+\|p^{ak}_0\|\not=0,
\end{equation}
\begin{equation}\label{e:dac1}
\al^k_ig_i(\ox^k_k-\ou^k_k)=0,\;\;i=1,\ldots,m,
\end{equation}
\begin{equation}\label{e:dac3}
\xi^{1k}_j(\|u^k_j\|-r_2-\ve_k)=0\;\;\mbox{for}\;\;j=0,\ldots,k,
\end{equation}
\begin{equation}
\label{e:dac3a}
\xi^{2k}_j(\|u^k_j\|-r_1+\ve_k)=0\;\;\mbox{for}\;\;j=0,\ldots,k,
\end{equation}
\begin{equation}\label{e:dac4}
\left\{\begin{array}{ll}
-p^{xk}_k\in\lm^k\partial\ph(\ox^k_k)-\disp\sum^m_{i=1}\al^k_i\nabla g_i(\ox^k_k-\ou^k_k),\\[1ex]
p^{uk}_k=\disp-\sum^m_{i=1}\al^k_i\nabla g_i(\ox^k_k-\ou^k_k)-2(\xi^{1k}_k+\xi^{2k}_k)\ou^k_k,\;\;p^{ak}_k=0,
\end{array}\right.
\end{equation}
\begin{equation}\label{e:dac5}
p^{uk}_{j+1}=\lm^k(v^{uk}_j+h^{-1}_k\th^{uk}_j),\;\; p^{ak}_{j+1}=\lm^k(v^{ak}_j+h^{-1}_k\th^{ak}_j),\;\;j=0,\ldots,k-1,
\end{equation}
\begin{equation}
\label{e:dac6}
\begin{aligned}
&\bigg(\dfrac{p^{xk}_{j+1}-p^{xk}_j}{h_k}-\lm^kw^{xk}_j,\dfrac{p^{uk}_{j+1}-p^{uk}_j}{h_k}-\lm^kw^{uk}_j,\dfrac{p^{ak}_{j+1}-p^{ak}_j}{h_k}-\lm^kw^{ak}_j,\\
&p^{xk}_{j+1}-\lm^k(v^{xk}_j+h^{-1}_k\th^{xk}_j)\bigg)\in \bigg(0,\dfrac{2}{h_k}(\xi^{1k}_j+\xi^{2k}_j)\ou^k_j,0,0\bigg)+N\bigg(\bigg(\ox^k_j,\ou^k_j,\oa^k_j,\dfrac{\ox^k_{j+1}-\ox^k_j}{-h_k}\bigg);\gph F\bigg)
\end{aligned}
\end{equation}
for $j=0,\ldots,k-1$ with the triples
\begin{equation}\label{e:dac9}
(\th^{xk}_j,\th^{uk}_j,\th^{ak}_j):=2\int^{t_{j+1}}_{t_j}\left(\dfrac{\ox^k_{j+1}-\ox^k_j}{h_k}-\dot{\ox}(t),\dfrac{\ou^k_{j+1}-\ou^k_j}{h_k}-\dot{\ou}(t),
\dfrac{\oa^k_{j+1}-\oa^k_j}{h_k}-\dot{\oa}(t)\right)dt
\end{equation}
and the subgradient collections
\begin{equation}
\label{e:dac7}
(w^{xk}_j,w^{uk}_j,w^{uk}_j,v^{xk}_j,v^{uk}_j,v^{ak}_j)\in\partial\ell_t\left(\oz^k_j,\dfrac{\oz^k_{j+1}-\oz^k_j}{h_k}\right),\;\;j=0,\ldots,k-1,
\end{equation}
where the sequence $\{\ve_k\}\dn0$ as $k\to\infty$ is taken from Theorem~{\rm\ref{Th:DA}}.
\end{theorem}\vspace*{-0.1in}
{\bf Proof.} For simplicity we drop the upper index ``$k$" in the notation below and consider the ``long" vector $y$ reflecting the collection of feasible solutions to each discrete-time problem $(P_k)$:
$$
y=(x_0,\ldots,x_k,u_0,\ldots,u_k,a_0,\ldots,a_k,X_0,\ldots,X_{k-1},U_0,\ldots,U_{k-1},A_0,\ldots,A_{k-1}).
$$
We now reduce $(P_k)$ to the following equivalent nondynamic problem of mathematical programming $(MP)$ with respect to the variable vector $y$, where the starting point $x_0$ is fixed:
\begin{eqnarray*}
\begin{aligned}
\mbox{minimize }\;&\varphi_{0}[y]:=\varphi(x_{k})+h_{k}\sum^{k-1}_{j=0}\ell(x_j,u_j,a_j,X_j,U_j,A_j)+\sum^{k-1}_{j=0}\int^{t_{j+1}}_{t_j}\left\|(X_j,U_j,A_j)-\dot{\oz}(t)
\right\|^{2}dt\\+&\left\|\frac{x^k_1-x^k_0}{h_k}-\dot\ox(0)\right\|^2+\dist^2\left(\left\|\dfrac{u^k_1-u^k_0}{h_k}\right\|;(-\infty,\Tilde{\mu}]\right)
+\dist^2\left(\sum^{k-2}_{j=0}\left\|U_{j+1}-U_j
\right\|;\big(-\infty,\Tilde{\mu}\big]\right)
\end{aligned}
\end{eqnarray*}
subject to the finitely many equality, inequality, and geometric constraints given by
\begin{eqnarray*}
\begin{aligned}
b^x_j(y):=&x_{j+1}-x_j-h_kX_j=0\;\mbox{ for }\;j=0,\ldots,k-1,\\
b^u_j(y):=&u_{j+1}-u_j-h_kU_j=0\;\mbox{ for }\;j=0,\ldots,k-1,\\
b^a_j(y):=&a_{j+1}-a_j-h_kA_j=0\;\mbox{ for }\;j=0,\ldots,k-1,\\
c_i(y):=&-g_i(x_k-u_k)\le 0\;\mbox{ for }\;i=1,\ldots,m,\\
d^1_j(y):=&\|u_j\|^2-(r_2+\ve_k)^2\le 0\;\mbox{ for }\;j=0,\ldots,k,\\
d^2_j(y):=&\|u_j\|^2-(r_1-\ve_k)^2\ge 0\;\mbox{ for }\;j=0,\ldots,k,\\
\phi_{j}(y):=&\left\|(x_j,u_j,a_j)-\oz(t_j)\right\|-\ve/2\le 0\;\mbox{ for }\;j=0,\ldots,k,\\
\phi_{k+1}(y):=&\sum^{k-1}_{j=0}\int^{t_{j+1}}_{t_j}\bigg(\left\|(X_j,U_j,A_j)-\dot{\oz}(t)\right\|^{2}\bigg)dt-\frac{\ve}{2}\le 0,\\
\phi_{k+2}(y):=&\sum^{k-2}_{j=0}\left\|U_{j+1}-U_j\right\|\le\Tilde{\mu}+1,\quad\phi_{k+3}(y):=\|u_1-u_0\|\le(\Tilde{\mu}+1)(t^k_1-t^k_0),\\
y\in\Xi_{j}:=&\big\{y|\;-X_{j}\in F(x_{j},u_{j},a_{j})\big\}\;\mbox{ for }\;j=0,\ldots,k-1,\\
y\in\Xi_{k}:=&\big\{y|\;x_{0}\;\mbox{ is fixed},\;(u_0,a_0)=\big(\ou(0),\oa(0)\big)\big\}.
\end{aligned}
\end{eqnarray*}
Let us apply the necessary optimality conditions from \cite[Theorem~5.24]{m-book2} to any local optimal solution $\oy$ of the finite-dimensional problem $(MP)$ written above with taking into account that by Theorem~\ref{Th:strong-DA-ilm} all the inequality constraints in $(MP)$ relating to functions $\phi_j$ as $j=0,\ldots,k+2$ are {\em inactive} for large $k\in\N$, and hence the corresponding multipliers do not appear in the optimality conditions. In this way we find dual elements $\lambda\ge 0,\;\alpha=(\al_1,\ldots,\al_m)\in\R^m_+,\;\xi^{1k}=(\xi^{1k}_0,\ldots,\xi^{1k}_k)\in\R^{k+1}_+,\xi^{2k}=(\xi^{2k}_0,\ldots,\xi^{2k}_k)\in\R^{k+1}_-,\;p_j=(p^x_j,p^u_j,p^a_j)\in\R^{2n+d}$ as $j=0,\ldots,k$, and
$$
y^*_j=(x^*_{0j},\ldots,x^*_{kj},u^*_{0j},\ldots,u^*_{kj},a^*_{0j},\ldots,a^*_{kj},X^*_{0j},\ldots,X^*_{(k-1)j},U^*_{0j},\ldots,U^*_{(k-1)j},A^*_{0j},\ldots,
A^*_{(k-1)j}),
$$
$j=0,\ldots,k$, which are not zero simultaneously, satisfy the conditions in \eqref{e:dac4} and the inclusions
\begin{equation}\label{e:dac10}
y^*_j\in N(\oy;\Xi_j)\;\;\mbox{for}\; j=0,\ldots,k,
\end{equation}
\begin{equation}\label{e:dac11}
-y^*_0-\ldots-y^*_k\in\lambda\partial\varphi_{0}(\oy)+\sum^m_{i=1}\alpha_i\nabla c_i(\oy)+\sum^{k}_{j=0}\xi^1_j\nabla d^1_j(\oy)+\sum^k_{j=0}\xi^2_j\nabla d^2_j(\oy)+\sum^{k-1}_{j=0}\nabla b_j(\oy)^*p_{j+1},
\end{equation}
\begin{equation*}
\al_i c_i(\oy)=0\;\mbox{ for }\;i=1,\ldots,m,
\end{equation*}
It is easy to see that the validity of \eqref{e:dac3} follows directly from the structure of the sets $\Xi_j$, and that the inclusions in \eqref{e:dac10} can be equivalently rewritten as
\begin{equation}\label{e:dac12}
\left(x^*_{jj},u^*_{jj},a^*_{jj},-X^{*}_{jj}\right)\in N\left(\left(\ox^{k}_{j},\ou^{k}_{j},\oa^{k}_{j},\dfrac{\ox^{k}_{j+1}-\ox^{k}_{j}}{-h_k}\right);\gph F\right),\;j=0,\ldots,k,
\end{equation}
while every other components of $y^*_j$ equals to zero. We conclude similarly that the only nonzero component of $y^*_k$ might be $(x^*_{0k},u^*_{0k},a^*_{0k})$. This gives us the equality
\begin{equation}\label{e:dac13}
\begin{aligned}
-y^*_0-y^*_1-\ldots-y^*_k=\big(&-x^*_{0k}-x^*_{00},-x^*_{11},\ldots,-x^*_{k-1,k-1},0,-u^*_{0k}-u^*_{00},\ldots,-u^*_{k-1,k-1},0,\\
&-a^*_{0k}-a^*_{00},-a^*_{11},\ldots,-a^*_{k-1,k-1},0,-X^*_{00},\ldots,-X^*_{k-1,k-1},0,\ldots,0\big).
\end{aligned}
\end{equation}
Next we calculate the sums on the right-hand side of \eqref{e:dac11}. It follows from the constructions above that
\begin{eqnarray*}
\begin{aligned}
\left(\sum^m_{i=1}\alpha_i\nabla c_i(\oy)\right)_{(x_k,u_k,a_k)}=&\left(-\sum^m_{i=1}\alpha_i\nabla g_i(x_k-u_k),\sum^m_{i=1}\alpha_i\nabla g_i(x_k-u_k),0\right),\\
\left(\sum^{k}_{j=0}\xi^1_j\nabla d^1_j(\oy)+\sum^{k}_{j=0}\xi^2_j\nabla d^2_j(\oy)\right)_{u_j}=&(2\xi^1_j+2\xi^2_j)\ou_j\;\mbox{ for }\;j=0,\ldots,k,\\
\left(\sum^{k-1}_{j=0}\big(\nabla b_j(\oy)\big)^*p_{j+1}\right)_{(x_j,u_j,a_j)}=&
\left\{\begin{array}{lcl}-p_1\quad\mbox{if }\;j=0,\\
p_j-p_{j+1}\quad\mbox{if }\;j=1,\ldots,k-1,\\
p_k\quad\mbox{if }\;j=k,
\end{array}\right.\\
\left(\sum^{k-1}_{j=0}\big(\nabla b_j(\oy)\big)^*p_{j+1}\right)_{(X,U,A)}=&-h_k p=\big(-h_kp^x_1,\ldots,-h_kp^x_k,-h_kp^u_1,\ldots,-h_kp^u_k,-h_kp^a_1,\ldots,-h_kp^a_k\big).
\end{aligned}
\end{eqnarray*}
Applying the subdifferential sum rule from \cite[Theorem~2.13]{m-book1} yields the inclusion
$$
\partial\varphi_0(\oy)\subset\partial\varphi(\ox_k)+h_k\sum^{k-1}_{j=0}\partial\ell_t(\ox_j,\ou_j,\oa_j,\bar{X}_j,
\bar{U}_j,\bar{A}_j)+\sum^{k-1}_{j=0}\nabla\rho_j(\oy)+\partial\sigma(\oy)
$$
with the real-valued functions $\rho_j(\cdot)$ and $\sigma(\cdot)$ given by
$$
\rho_j(y):=\int^{t_{j+1}}_{t_j}\left\|(X_j,U_j,A_j)-\dot{\oz}(t)\right\|^2dt,
$$
$$
\sigma(y):=\dist^2\left(\left\|\dfrac{u^k_1-u^k_0}{h_k}\right\|;\big(-\infty,\Tilde{\mu}\big]\right)+\dist^2\left(\sum^{k-2}_{j=0}\left\|U_{j+1}-U_j
\right\|;\big(-\infty,\Tilde{\mu}\big]\right).
$$
Using now the differentiability of $\psi(x):=\dist^2(x;(-\infty,\Tilde{\mu}])$ with the gradient $\nabla\psi(x)=0$ whenever $x\le\Tilde{\mu}$ and combining it with second condition in \eqref{e:discrete-estimate} tells us that $\partial\sigma(\oy)=\{0\}$. Furthermore, we get
$$
\nabla\rho_j(\oy)=\nabla_{X_j,U_j,A_j}\rho(\oy)=(\theta^x_j,\theta^u_j,\theta^a_j)
$$
for nonzero components. Putting this together shows that $\lambda\partial\varphi_0(\oy)$ in \eqref{e:dac11} is represented as
$$
\begin{array}{ll}\lambda(h_kw^x_0,h_kw^k_1,\ldots,h_kw^x_{k-1},\vartheta^k,h_kw^u_0,\ldots,h_kw^u_{k-1},0,h_kw^a_0,\ldots,h_kw^a_{k-1},0,\theta^x_0+h_kv^x_0,\ldots,\\
\theta^x_{k-1}+h_kv^x_{k-1},\theta^u_0+h_kv^u_0,\ldots,\theta^u_{k-1}+h_kv^u_{k-1},\theta^a_0+h_kv^a_0,\ldots,\theta^a_{k-1}+h_kv^a_{k-1})
\end{array}
$$
with $\vartheta^k\in\partial\varphi(\ox_k)$ and with the components of $(w^x,w^u,w^a,v^x,v^u,v^a)$ satisfying \eqref{e:dac7}. Involving \eqref{e:dac13}, we derive from \eqref{e:dac10} the following relationships:
\begin{equation*}
\begin{aligned}
-x^*_{0k}-x^*_{00}=&\lambda h_kw^x_0-p^x_1,\\
-x^*_{jj}=&\lambda h_kw^x_j+p^x_j-p^x_{j+1}\;\mbox{ for }\;j=1,\ldots,k-1,\\
0=&\lambda\vartheta^k-\sum^m_{i=1}\alpha_i\nabla_x\Tilde g_i(x_k,u_k)+p^x_k,\\
-u^*_{0k}-u^*_{00}=&\lambda h_kw^u_0+(2\xi^1_0+2\xi^2_0)\ou_0-p^u_1,\\
-u^*_{jj}=&\lambda h_kw^u_j+2(\xi^1_j+\xi^2_j)\ou_j+p^u_j-p^u_{j+1}\;\mbox{ for }\;j=1,\ldots,k-1,\\
0=&\sum^m_{i=1}\alpha_i\nabla_u\Tilde g_i(x_k,u_k)+p^u_k+(2\xi^1_k+2\xi^2_k)\ou_k,\\
-a^*_{0k}-a^*_{00}=&\lambda h_kw^a_0-p^a_1,\\
-a^*_{jj}=&\lambda h_kw^a_j+p^a_j-p^a_{j+1}\;\mbox{ for }\;j=1,\ldots,k-1,\\
0=&p^a_k,\\
-X^*_{jj}=&\lambda(h_kv^x_j+\theta^x_j)-h_kp^x_{j+1}\;\mbox{ for }\;j=0,\ldots,k-1,\\
0=&\lambda(h_kv^u_j+\theta^u_j)-h_kp^u_{j+1}\;\mbox{ for }\;j=0,\ldots,k-1,\\
0=&\lambda(h_kv^a_j+\theta^a_j)-h_kp^a_{j+1}\;\mbox{ for }\;j=0,\ldots,k-1.
\end{aligned}
\end{equation*}
\vspace*{-0.05in}
To complete the proof, we proceed similarly to the last part in the proof of \cite[Theorem~7.1]{cm1}.$\h$

Employing the second-order calculations conducted in Section~6 allows us to derive from Theorem~\ref{Th:EL-DA} necessary optimality conditions for each discrete control problem $(P_k)$ expressed entirely via its given data. For simplicity of the formulation we assume that the perturbation mapping $f$ is smooth with respect to the state variable $x$. The reader can certainly proceed with the case of Lipschitz continuous mapping $f$ by using well-developed calculus rules for our basic first-order generalized differential constructions \cite{m-book1}.\vspace*{-0.1in}

\begin{theorem}{\bf (optimality conditions for discretized sweeping control problems via their original data).}\label{Th:OC-DP} Given an optimal control $\oz^k=(\ox^k,\ou^k,\oa^k)$ to discrete-time problem $(P_k)$ with any fixed $k\in\N$ and with the sweeping velocity mapping $F$ defined in \eqref{e:vel-map},
suppose that the functions $g_i$ in \eqref{e:mset} are of class ${\cal C}^2$ and the perturbation mapping $f(\cdot,a)$ is of class ${\cal C}^1$ around the optimal points. Then there are dual elements $(\lm^k,\be^k,\xi^{1k},\xi^{2k},p^k)$ together with vectors $\eta^k_j\in\R^m_+$ as $j=0,\ldots,k$ and $\gg^k_j\in\R^m$ as $j=0,\ldots,k-1$ satisfying \eqref{e:dac3} and such that the following conditions hold:\\[1ex]
{\sc nontriviality condition}
\begin{equation}\label{e:dac14}
\lm^k+\|\eta^k_k\|+\|\xi^{1k}+\xi^{2k}\|+\sum^{k-1}_{j=0}\|p^{xk}_j\|+\|p^{uk}_0\|+\|p^{ak}_0\|\not=0;
\end{equation}
{\sc primal-dual dynamic relationship} for all $j=0,\ldots,k-1$:
\begin{equation}\label{e:dac15}
\dfrac{\ox^k_{j+1}-\ox^k_j}{h_k}+f(\ox^k_j,\oa^k_j)=\sum_{i\in I(\ox^k_j-\ou^k_j)}\eta^k_{ji}\nabla g_i(\ox^k_j-\ou^k_j),
\end{equation}
\begin{equation}\label{e:dac16}
\begin{aligned}
&\dfrac{p^{xk}_{j+1}-p^{xk}_j}{h_k}-\lm^kw^{xk}_j=\nabla_xf(\ox^k_j,\oa^k_j)^*(\lm^k(v^{xk}_j+h^{-1}_k\th^{xk}_j)-p^{xk}_{j+1})\\
&-\sum^m_{i=1}\eta^k_{ji}\nabla^2g_i(\ox^k_j-\ou^k_j)(\lm^k(v^{xk}_j+h^{-1}_k\th^{xk}_j)-p^{xk}_{j+1})-\sum^m_{i=1}\gg^k_{ji}\nabla g_i(\ox^k_j-\ou^k_j),
\end{aligned}
\end{equation}
\begin{equation}\label{e:dac17}
\begin{array}{ll}
\dfrac{p^{uk}_{j+1}-p^{uk}_j}{h_k}-\lm^kw^{uk}_j-\dfrac{2}{h_k}(\xi^{1k}_j+\xi^{2k}_j)\ou^k_j&=\disp\sum^m_{i=1}\eta^k_{ji}\nabla^2g_i(\ox^k_j
-\ou^k_j)(\lm^k(v^{xk}_j+h^{-1}_k\th^{xk}_j)-p^{xk}_{j+1})\\
&+\disp\sum^m_{i=1}\gg^k_{ji}\nabla g_i(\ox^k_j-\ou^k_j),
\end{array}
\end{equation}
\begin{equation}
\label{e:dac18}
\begin{aligned}
&\dfrac{p^{ak}_{j+1}-p^{ak}_j}{h_k}-\lm^kw^{ak}_j=\nabla_af(\ox^k_j,\oa^k_j)^*(\lm^k(v^{xk}_j+h^{-1}_k\th^{xk}_j)-p^{xk}_{j+1}),
\end{aligned}
\end{equation}
where $(w^{xk}_j,w^{uk}_j,w^{ak}_j,v^{xk}_j,v^{uk}_j,v^{ak}_j)$ are taken from \eqref{e:dac7} while the active constraint index set $I(\cdot)$ and the triples $\th^{xk}_j,\th^{uk}_j,\th^{ak}_j$ are defined in \eqref{e:a-index} and \eqref{e:dac9}, respectively;\\[1ex]
{\sc transversality conditions}
\begin{equation}\label{e:dac19}
\left\{\begin{array}{ll}
-p^{xk}_k\in\lm^k\partial\ph(\ox^k_k)-\disp\sum^m_{i=1}\eta^k_{ki}\nabla g_i(\ox^k_k-\ou^k_k),\\[1ex]
p^{uk}_k=\disp-\sum^m_{i=1}\eta^k_{ki}\nabla g_i(\ox^k_k-\ou^k_k)-2(\xi^{1k}_k+\xi^{2k}_k)\ou^k_k,\;\; p^{ak}_k=0
\end{array}\right.
\end{equation}
with dual vectors $\xi^{1k}_k$ and $\xi^{2k}_k$ satisfying the inclusions
\begin{equation}\label{e:dac22a}
\xi^{1k}_k\in N_{[0,r_2+\ve_k]}(\|\ou^k_k\|),\;\;\xi^{2k}_k\in N_{[r_1-\ve_k,\infty)}(\|\ou^k_k\|);
\end{equation}
{\sc complementarity slackness conditions}
\begin{equation}
\label{e:dac20}
[g_i(\ox^k_j-\ou^k_j)>0]\Lto\eta^k_{ji}=0,
\end{equation}
\begin{equation}
\label{e:dac21}
\left\{
\begin{array}{ll}
[i\in I_1(-p^{xk}_{j+1}+\lm^k(h^{-1}_k\th^{xk}_j+v^{xk}_j))],\; \mbox{i.e.,}\; [g_i(\ox^k_j-\ou^k_j)>0\;\mbox{or} \\
\eta^k_{ji}=0, \la \nabla g_i(\ox^k_j-\ou^k_j),-p^{xk}_{j+1}+\lm^k(h^{-1}_k\th^{xk}_j+v^{xk}_j)\ra>0]\Lto[\gg^k_{ji}=0],
\end{array}\right.
\end{equation}
\begin{equation}
\label{e:dac22}
\left\{
\begin{array}{ll}
[i\in I_2(-p^{xk}_{j+1}+\lm^k(h^{-1}_k\th^{xk}_j+v^{xk}_j))],\; \mbox{i.e.,}\; [g_i(\ox^k_j-\ou^k_j)=0,\eta^k_{ji}=0,\;\mbox{and}\\
\la \nabla g_i(\ox^k_j-\ou^k_j),-p^{xk}_{j+1}+\lm^k(h^{-1}_k\th^{xk}_j+v^{xk}_j)\ra<0]\Lto[\gg^k_{ji}\ge0]
\end{array}\right.
\end{equation}
for $j=0,\ldots,k-1$ and $i=1,\ldots,m$, where the index subsets $I_1(\cdot)$ and $I_2(\cdot)$ are taken from \eqref{e:indexes}, in addition to
\eqref{e:dac3} together with the implications
\begin{equation}
\label{e:dac23}
[g_i(\ox^k_j-\ou^k_j)>0]\Lto\gg^k_{ji}=0\;\mbox{for}\;\;j=0,\ldots,k-1\;\;\mbox{and}\;\;i=1,\ldots,m,
\end{equation}
\begin{equation}\label{e:dac24}
[g_i(\ox^k_k-\ou^k_k)>0]\Lto\eta^k_{ki}=0\;\;\mbox{for}\;\;i=1,\ldots,m,\;\mbox{ and}
\end{equation}
\begin{equation}\label{e:dac25}
\eta^k_{ji}>0\Lto[\la\nabla g_i(\ox^k_j-\ou^k_j),-p^{xk}_{j+1}+\lm^k(h^{-1}_k\th^{xk}_j+v^{xk}_j)\ra=0].
\end{equation}
Furthermore, assuming the surjectivity of the Jacobian matrix $\{\nabla g(\ox^k_j-\ou^k_j)\}$ ensures the validity of the {\sc enhanced nontriviality condition}
\begin{equation}
\label{e:dac26}
\lm^k+\|\xi^{1k}+\xi^{2k}\|+\|p^{uk}_0\|\not=0.
\end{equation}
\end{theorem}\vspace*{-0.1in}
{\bf Proof.} It follows from \eqref{e:dac6} and the coderivative definition \eqref{e:cor} that
\begin{equation}\label{e:dac27}
\begin{aligned}
&\bigg(\dfrac{p^{xk}_{j+1}-p^{xk}_j}{h_k}-\lm^kw^{xk}_j,\dfrac{p^{uk}_{j+1}-p^{uk}_j}{h_k}-\lm^kw^{uk}_j-\dfrac{2}{h_k}(\xi^{1k}_j+\xi^{2k}_j)\ou^k_j,
\dfrac{p^{ak}_{j+1}-p^{ak}_j}{h_k}-\lm^kw^{ak}_j\bigg)\\
&\in D^*F\bigg(\ox^k_j,\ou^k_j,\oa^k_j,\dfrac{\ox^k_{j+1}-\ox^k_j}{-h_k}\bigg)(\lm^k(h^{-1}_k\th^{xk}_j+v^{xk}_j)-p^{xk}_{j+1}),\;\;j=0,\ldots,k-1.
\end{aligned}
\end{equation}
By $\dfrac{\ox^k_{j+1}-\ox^k_j}{-h_k}-f(\ox^k_j,\oa^k_j)\in N(\ox^k_j-\ou^k_j;C)$ for $j=0,\ldots,k-1$ and the representation of $F$ in \eqref{e:F-rep}, we find  vectors $\eta^k_j\in\R^m_+$, $j=0,\ldots,k-1$, such that the conditions in \eqref{e:dac15} and \eqref{e:dac20} are satisfied. Employing the second-order upper estimate from Theorem~\ref{Th:co-cal} with $x:=\ox^k_j,u:=\ou^k_j,a:=\oa^k_j,w:=\dfrac{\ox^k_{j+1}-\ox^k_j}{-h_k}$, and $y:=\lm^k(h^{-1}_k\th^{xk}_j+v^{xk}_j)-p^{xk}_{j+1}$ and combining it with \eqref{e:dom-co-est} give us $\gg^k_j\in\R^m$ for which
$$
\begin{aligned}
&\bigg(\dfrac{p^{xk}_{j+1}-p^{xk}_j}{h_k}-\lm^kw^{xk}_j,\dfrac{p^{uk}_{j+1}-p^{uk}_j}{h_k}-\lm^kw^{uk}_j-\dfrac{2}{h_k}(\xi^{1k}_j+\xi^{2k}_k)\ou^k_j,\dfrac{p^{ak}_{j+1}-p^{ak}_j}{h_k}-\lm^kw^{ak}_j\bigg)\\
=&\bigg(\nabla_xf(\ox^k_j,\oa^k_j)^*(\lm^k(v^{xk}_j+h^{-1}_k\th^{xk}_j)-p^{xk}_{j+1})-\sum^m_{i=1}\eta^k_{ji}\nabla^2g_i(\ox^k_j-\ou^k_j)(\lm^k(v^{xk}_j+h^{-1}_k\th^{xk}_j)-p^{xk}_{j+1})\\
-&\sum^m_{i=1}\gg^k_{ji}\nabla g_i(\ox^k_j-\ou^k_j),\sum^m_{i=1}\eta^k_{ji}\nabla^2g_i(\ox^k_j-\ou^k_j)(\lm^k(v^{xk}_j+h^{-1}_k\th^{xk}_j)-p^{xk}_{j+1})+\sum^m_{i=1}\gg^k_{ji}\nabla g_i(\ox^k_j-\ou^k_j),\\
&\nabla_af(\ox^k_j,\oa^k_j)^*(\lm^k(v^{xk}_j+h^{-1}_k\th^{xk}_j)-p^{xk}_{j+1})\bigg),\quad j=0,\ldots,k-1.
\end{aligned}
$$
This yields the validity of all the conditions in \eqref{e:dac16}, \eqref{e:dac17}, \eqref{e:dac18}, \eqref{e:dac21}, and \eqref{e:dac22}. Put $\eta^k_k:=\al_k$ with $\al_k$ taken from the statement of Theorem~\ref{Th:EL-DA} and observe that $\eta^k_j\in\R^m_+$ for $j=0,\ldots,k$. In this way we deduce the nontriviality condition \eqref{e:dac14} from \eqref{e:nontri-da} and the transversality conditions \eqref{e:dac19} from \eqref{e:dac4}. Furthermore, \eqref{e:dac24} follows immediately from \eqref{e:dac1} and the definition of $\eta^k_k$, while \eqref{e:dac27} yields
$$
\lm^k(h^{-1}_k\th^{xk}_j+v^{xk}_j)-p^{xk}_{j+1}\in\dom D^*N_C\bigg(\ox^k_j-\ou^k_j,\dfrac{\ox^k_{j+1}-\ox^k_j}{-h_k}-f(\ox^k_j,\oa^k_j)\bigg),
$$
This implies by \eqref{e:dom-co-est} that \eqref{e:dac25} holds. Inclusions \eqref{e:dac22a} follow from \eqref{e:dac3} and \eqref{e:dac3a}.\vspace*{-0.05in}

To complete the proof of the theorem, it remains to justify the enhanced nontriviality condition \eqref{e:dac26} under the surjectivity of the Jacobians $\{\nabla g(\ox^k_j-\ou^k_j)\}$. Arguing by contradiction, suppose that \eqref{e:dac26} fails, i.e., $\lm^k=0,\;\xi^{1k}+\xi^{2k}=0$, and $p^{uk}_0=0$. Then it follows from \eqref{e:dac5} that $p^{uk}_j=0$ for $j=0,\ldots,k$ and $p^{ak}_j=0$ for $j=1,\ldots,k$. Employing the second condition in \eqref{e:dac19} with $p^{uk}_k=0$ tells us that $\disp\sum^m_{i=1}\eta^k_{ki}\nabla g_i(\ox^k_k-\ou^k_k)=0$, and so $p^{xk}_k=0$ by the first condition therein. We also get from \eqref{e:dac17} that
$$
\sum^m_{i=1}\eta^k_{ji}\nabla^2g_i(\ox^k_j-\ou^k_j)(\lm^k(v^{xk}_j+h^{-1}_k\th^{xk}_j)-p^{xk}_{j+1})+\sum^m_{i=1}\gg^k_{ji}\nabla g_i(\ox^k_j-\ou^k_j)=0,\quad j=0,\ldots,k-1.
$$
Using this together with \eqref{e:dac16} and $p^{xk}_k=0$ shows that $p^{xk}_j=0$ for all $j=0,\ldots,k-1$. Finally, it follows from \eqref{e:dac18} that $p^{ak}_0=0$, which contradicts the validity of \eqref{e:dac14} and thus verifies \eqref{e:dac26}. $\h$\vspace*{-0.2in}

\section{Necessary Conditions for Sweeping Optimal Solutions}
\setcounter{equation}{0}\vspace*{-0.1in}

This section is the culmination of the paper. Given an arbitrary {\em relaxed intermediate local minimizer} $\oz(\cdot)$ of the sweeping optimal control problem $(P)$ and using the discrete approximation method, we derive verifiable necessary optimality conditions for $\oz(\cdot)$, which expressed entirely via the problem data, by combining the strong convergence of discrete optimal solutions to $\oz(\cdot)$ established Section~5 and the necessary optimality conditions in discrete approximations taken from Section~7 with a rather involved technique to justify a proper convergence of the adjoint trajectories from Theorem~\ref{Th:OC-DP}. The latter technique developed here is heavily based on the underlying properties of our basic generalized differential constructions and the second-order calculations of Section~6. As discussed in Section~4, {\em no relaxation} is needed if $\oz(\cdot)$ is a {\em strong} local minimizer of $(P)$ and either controls are located only in perturbations, or the set $C$ in \eqref{e:mset} is convex in addition to the standing assumptions formulated in Section~2.\vspace*{-0.05in}

We now add to the standing assumptions the one on time dependence of the basic subdifferential $\partial\ell$ of the running cost in \eqref{e:Bolza} taken below with respect to all but $t$ variables. It is well know that the subdifferential mapping \eqref{e:sub} is {\em robust} (which reduces to the graph-closedness for continuous functions) with respect to the variables of subdifferentiation. We suppose that this robustness property keeps holding when the time parameter is involved into the limiting procedure. Precisely it amounts to saying that
$$
\partial\ell\big(t,\oz(t),\dot{\oz}(t)\big)=\disp\Limsup_{(\tau,u,v)\st{\ell}{\to}(t,z(t),\dot{z}(t))}\partial\ell(\tau,u,v)\;\mbox{ a.e. }\;t\in[0,T]
$$
around the given local optimal solution to $(P)$, where ``$\Limsup$" stands for the Painlev\'e-Kuratowski outer/upper limit \cite{rw}. This assumption is not restrictive and is satisfied, in particular, for smooth functions with time-continuous derivatives as well as in broad nonsmooth settings; see \cite{m95,m-book2}.\vspace*{-0.1in}

\begin{theorem}{\bf (optimality conditions for the nonconvex sweeping process.)}\label{Th:NOC} Given an r.i.l.m.\ $\oz(\cdot)=(\ox(\cdot),\ou(\cdot),\oa(\cdot))$ for problem $(P)$, suppose in addition to the standing assumptions and those in Theorem~{\rm\ref{Th:DA}} that $\ell$ is continuous in $t$ a.e.\ on $[0,T]$ and is represented as
\begin{equation}\label{separ}
\ell(t,z,\dot z)=\ell_1(t,z,\dot x)+\ell_2(t,\dot u)+\ell_3(t,\dot a)
\end{equation}
where the local Lipschitz constants of $\ell_1(t,\cdot,\cdot)$ and $\ell_3(t,\cdot)$ are essentially bounded on $[0,T]$ and continuous at a.e.\ $t\in[0,T]$ including $t=0$, and where $\ell_2$ is differentiable in $\dot u$ on $\R^n$ with the estimates
\begin{equation}\label{estimate1}
\begin{array}{c}
\|\nabla_{\dot{u}}\ell_2(t,\dot{u},\dot{a})\|\le L\|\dot{u}\|\;\mbox{ and }\;\|\nabla_{\dot{u}}\ell_2(t,\dot{u}_1)-\nabla_{\dot{u}}\ell_2(s,\dot{u}_2)\|\le L|t-s|+L\|\dot{u}_1-\dot{u}_2\|
\end{array}
\end{equation}
holding for all $t,s\in[0,T]$, $\dot a\in\R^d$, and $\dot u,\dot u_1,\dot u_2\in\R^n$ with some uniform constant $L>0$. Then there are $\lm\ge 0$, $p(\cdot)=(p^x(\cdot),p^u(\cdot),p^a(\cdot))\in W^{1,2}([0,T];\R^n\times\R^n\times\R^d)$, $w(\cdot)=(w^x(\cdot),w^u(\cdot),w^a(\cdot))\in L^2([0,T];\R^{2n+d})$, and $v(\cdot)=(v^x(\cdot),v^u(\cdot),v^a(\cdot))\in L^2([0,T];\R^{2n+d})$ well defined at $t=0$ and satisfying
\begin{equation}\label{subg}
\big(w(t),v(t)\big)\in\co\partial\ell\big(t,\oz(t),\dot{\oz}(t)\big)\;\mbox{ a.e. }\;t\in[0,T]
\end{equation}
as well as measures $\gamma=(\gg_1,\ldots,\gg_n)\in C^*([0,T];\R^n)$, $\xi^1\in C^*([0,T];\R_+)$, and $\xi^2\in C^*([0,T];\R_-)$ on $[0,T]$ such that the following conditions hold:\\
$\bullet$ {\sc Primal-dual dynamic relationships:}
\begin{equation}\label{e:pddr}
\dot{\ox}(t)+f\big(\ox(t),\oa(t)\big)=\sum_{i=1}^m\eta_i(t)\nabla g_i\big(\ox(t)-\ou(t)\big)\;\mbox{ for a.e. }\;t\in[0,T]
\end{equation}
with $\eta(\cdot)\in L^2([0,T];\R^+)$ a.e.\ uniquely determined by representation \eqref{e:pddr} and well defined at $t=T$;
\begin{equation}
\label{e:EL}\dot p(t)=\lm w(t)+\left(\nabla_xf(\ox(t),\oa(t))^*(\lm v^x(t)-q^x(t)),0,\nabla_af(\ox(t),\oa(t))^*(\lm v^x(t)-q^x(t))\right),
\end{equation}
\begin{equation}\label{e:qu-qa}
q^u(t)=\lm\nabla_{\dot u}\ell\big(t,\dot{\ou}(t)\big),\;\;q^a(t)\in\lm\partial_{\dot a}\ell_3\big(t,\dot{\oa}(t)\big)\;\mbox{ a.e. }\;t\in[0,T],
\end{equation}
where $q=(q^x,q^u,q^a):[0,T]\to\R^n\times\R^n\times\R^d$ is a vector function of bounded variation, and its left-continuous representative is given for all $t\in[0,T]$, except at most a countable subset, by
\begin{equation}\label{e:p-q relation}
q(t):=p(t)-\int_{[t,T]}\left(-d\gg(s),2\ou(s)d(\xi^1(s)+\xi^2(s))+d\gg(s),0\right).
\end{equation}
Furthermore, for a.e.\ $t\in[0,T]$ including $t=T$ and for all $i=1,\ldots,m$ we have
\begin{equation}\label{e:implications}
g_i\big(\ox(t)-\ou(t)\big)>0\Lto\eta_i(t)=0,\;\;\eta_i(t)>0\Lto\la\nabla g_i\big(\ox(t)-\ou(t),\lm v^x(t)-q^x(t))\big\ra=0.
\end{equation}
$\bullet$ {\sc Transversality conditions}
\begin{equation}\label{e:right-trans}
\left\{\begin{array}{ll}
-p^x(T)+\disp\sum_{i\in I(\ox(T)-\ou(T))}\eta_i(T)\nabla g_i\big(\ox(T)-\ou(T)\big)\in\lm\partial\ph\big(\ox(T)\big),\quad p^a(T)=0,\\[1ex]
p^u(T)\disp-\sum_{i\in I(\ox(T)-\ou(T))}\eta_i(T)\nabla g_i(\ox(T)-\ou(T))\in
-2\ou(T)\left(N_{[0,r_2]}(\|\ou(T)\|)+N_{[r_1,\infty)}(\|\ou(T)\|)\right)\\[1ex]
\end{array}\right.
\end{equation}
with the validity of the inclusion
\begin{equation}\label{e:right-trans1}
-\sum_{i\in I(\ox(T)-\ou(T))}\eta_i(T)\nabla g_i\big(\ox(T)-\ou(T)\big)\in N_C\big(\ox(T)-\ou(T)\big).
\end{equation}
$\bullet$ {\sc Measure nonatomicity conditions:}\\
{\bf(a)} Take $t\in[0,T]$ with $g_i(\ox(t)-\ou(t))>0$ whenever $i=1,\ldots,m$. Then there is a neighborhood $V_t$ of $t$ in $[0,T]$ such that $\gg(V)=0$ for all the Borel subsets $V$ of $V_t$.\\
{\bf(b)} Take $t\in[0,T]$ with $r_1<\|\ou(t)\|<r_2$. Then there is a neighborhood $W_t$ of $t$ in $[0,T]$ such that $\xi^1(W)=0$ and $\xi^2(W)=0$ for all  the Borel subsets $W$ of $W_t$.\\
$\bullet$ {\sc Nontriviality conditions:} We always have:
\begin{equation}\label{e:nontri-condition}
\lm+\|q^u(0)\|+\|p(T)\|+\|\xi^1\|_{TV}+\|\xi^2\|_{TV}>0.
\end{equation}
Furthermore, the following implications hold while ensuring the {\sc enhanced nontriviality}:
\begin{equation}\label{e:enhanced-nontri}
\big[g_i\big(x_0-\ou(0)\big)>0,\;i=1,\ldots,m\big]\Longrightarrow\big[\lm+\|p(T)\|+\|\xi^1\|_{TV}+\|\xi^2\|_{TV}>0\big],
\end{equation}
\begin{equation}\label{e:enhanced-nontri1}
\big[g_i\big(\ox(T)-\ou(T)\big)>0,\;r_1<\|\ou(T)\|<r_2,\;i=1,\ldots,m\big]\Longrightarrow\big[\lm+\|q^u(0)\|+\|\xi^1\|_{TV}+\|\xi^2\|_{TV}>0\big],
\end{equation}
where $\|\xi\|_{TV}$ stands for the measure total variation on $[0,T]$.
\end{theorem}\vspace*{-0.1in}
{\bf Proof.} We split the proof into the following major steps in accordance to the statement of the theorem.\\\vspace*{0.02in}
{\bf Step~1:} {\em Subgradients of the running cost.} To verify the subdifferential inclusion \eqref{subg}, take the subgradient sequence $\{w^k_j,v^k_j\}$ from Theorem~\ref{Th:OC-DP} and consider the piecewise constant extensions $w^k,v^k\colon[0,T]\to\R^{2n+d}$. It follows from \eqref{e:dac7} therein. The imposed assumptions and the structure of $\ell$ in \eqref{separ} with estimates \eqref{estimate1} ensure that the subgradient sets $\partial\ell(t,\cdot)$ are uniformly $L^2$-bounded near $\oz(\cdot)$, and hence the sequence $\{(w^k(\cdot),v^k(\cdot))\}$ is weakly compact in $L^2([0,T];\R^{2(2n+d)}):=L^2[0,T]$. Without relabeling we get the weak convergence $(w^k(\cdot),v^k(\cdot))\to(w(\cdot),v(\cdot))\in L^2[0,T]$ and thus, by Mazur's theorem, the strong $L^2$-convergence to $(w(\cdot),v(\cdot))$ of a sequence of convex combinations of $(w^k(\cdot),v^k(\cdot))$; that is, the a.e.\ convergence on $[0,T]$ to the above limiting pair of some subsequence of the latter. This readily verifies \eqref{subg} by taking into account the assumed a.e.\ continuity of $\ell$ in $t$ and the robustness of its subdifferential.\\\vspace*{0.02in}
{\bf Step~2:} {\em Verification of the primal dynamic limiting relationships.} We claim two of them: the differential equation \eqref{e:pddr} and the first implication in \eqref{e:implications}. We proceed by passing to the limit in \eqref{e:dac15}, \eqref{e:dac25} and first construct the piecewise constant functions on $[0,T]$ by
$$
\theta^{k}(t):=\dfrac{\theta^{k}_j}{h_k}\;\mbox{ as }\;t\in[t^k_j,t^k_{j+1}),\;j=0,\ldots,k-1,\quad k\in\N,
$$
where $\theta^{k}_j$ are taken from \eqref{e:dac9}. It follows from Theorem~\ref{Th:OC-DP} that
\begin{equation}\label{theta}
\begin{aligned}
\int^T_0\|\theta^{xk}(t)\|^2dt=&\sum^{k-1}_{j=0}\dfrac{\|\theta^{xk}_j\|^2}{h_k}\le\dfrac{4}{h_k}\sum^{k-1}_{j=0}\bigg(\int^{t^k_{j+1}}_{t^k_j}\Big\|\dot{\ox}(t)
-\dfrac{\ox^k_{j+1}-\ox^k_j}{h_k}\Big\|dt\bigg)^2\\\le&4\sum^{k-1}_{j=0}\bigg(\int^{t^k_{j+1}}_{t^k_j}\left\|\dot{\ox}(t)-\dfrac{\ox^k_{j+1}
-\ox^k_j}{h_k}\right\|^2\bigg)dt=4\int^T_0\left\|\dot{\ox}(t)-\dot{\ox}^k(t)\right\|^2dt\to 0
\end{aligned}
\end{equation}
as $k\to\infty$ with the same conclusion for $\theta^{uk}(\cdot)$ and $\theta^{ak}(\cdot)$. Thus some subsequences of these functions converge to zero a.e.\ on $[0,T]$. Invoking further the vectors $\eta^k_j\in\R^m_+$ from Theorem~\ref{Th:OC-DP}, define the piecewise constant functions $\eta^k(\cdot)$ on $[0,T]$ by $\eta^k(t):=\eta^k_j$ as $t\in[t^k_j,t^k_{j+1})$ with $\eta^k(T):=\eta^k_k$ and deduce from \eqref{e:dac15} for each $k\in\N$ we have the relationships
\begin{equation}\label{e:dpddr}
\dot\ox^k(t)+f\big(\ox^k(t),\oa^k(t)\big)=\sum^m_{i=1}\eta^k_i(t)\nabla g_i\big(\ox^k(t)-\ou^k(t)\big)\;\mbox{ if }\;t\in(t^k_j,t^k_{j+1}).
\end{equation}
The feasibility of $\oz(\cdot)$ in $(P)$ tells us that $-\dot\ox(t)\in N_C(\ox(t)-\ou(t))+f(\ox(t),\oa(t))$ for a.e.\ $t\in[0,T]$, where the closed-valued normal cone mapping $N_C(\cdot)$ is measurable by \cite[Theorem~14.26]{rw}. The classical measurable selection result (see, e.g., \cite[Corollary~14.6]{rw}) gives us nonnegative measurable functions $\eta_i(\cdot)$ on $[0,T]$ as $i=1,\ldots,m$ for which the differential equation \eqref{e:pddr} and the first implication in \eqref{e:implications} are satisfied. Then invoking \eqref{e:dpddr} and \eqref{e:pddr} yields the equalities
$$
\dot\ox(t)-\dot\ox^k(t)+f\big(\ox(t),\oa(t)\big)-f\big(\ox^k(t),\oa^k(t)\big)=\sum^m_{i=1}\left[\eta_i(t)\nabla g_i\big(\ox(t)-\ou(t)\big)-\eta^k_i(t)\nabla g_i\big(\ox^k(t)-\ou^k(t)\big)\right]
$$
whenever $t\in(t^k_j,t^k_{j+1})$ and $j=0,\ldots,k-1$, which imply the estimate
\begin{equation*}
\Big\|\sum^m_{i=1}\big[\eta_i(t)\nabla g_i\big(\ox(t)-\ou(t)\big)-\eta^k_i(t)\nabla g_i\big(\ox^k(t)-\ou^k(t)\big)\big]\Big\|\le\big\|\dot\ox(t)-\dot\ox^k(t)\big\|+\big\|f\big(\ox(t),\oa(t)\big)-f\big(\ox^k(t),\oa^k(t)\big)\big\|
\end{equation*}
on $(t^k_j,t^k_{j+1})$. Passing to the limit as $k\to\infty$ in this estimate with replacing $t^k_j$ by
\begin{equation}\label{max-time}
\nu^k(t):=\max\big\{t^k_j\big|\;t^k_j\le t,\;0\le j\le k\big\},\quad t\in[0,T],
\end{equation}
and taking into account that $I(\ox^k(\cdot)-\ou^k(\cdot))\subset I(\ox(\cdot)-\ou(\cdot))$ for $k\in\N$ sufficiently  large and that the sequence $\{\oz^k(\cdot)\}$ converges to $\ox(\cdot)$ strongly in $W^{1,2}([0,T])$, we get
$$
\sum_{i\in I(\ox(t)-\ou(t))}\left[\eta_i(t)\nabla g_i\big(\ox(t)-\ou(t)\big)-\eta^k_i(t)\nabla g_i\big(\ox^k(t)-\ou^k(t)\big)\right]\to 0\;\mbox{ a.e. }\;t\in[0,T].
$$
Observe also that $\eta^k(\cdot)\to\eta(\cdot)$ a.e.\ on $[0,T]$ by \eqref{e:bound-grad} and \eqref{e:w-inverse-triangle-in}. Postponing till Step~5 the verification of the claim that the sequence $\{\eta^k_k\}$ converges to the well-defined vector $(\eta_1(T),\ldots,\eta_m(T))$, we check now that $\eta(\cdot)\in L^2([0,T];\R^m_+)$. In fact, it follows from the estimates
$$
\begin{aligned}
\eta_i(t)\le\dfrac{1}{M_1}\eta_i(t)\big\|\nabla g_i\big(\ox(t)-\ou(t)\big)\big\|&\le\dfrac{1}{M_1}\sum_{i\in I(\ox(t)-\ou(t))}\eta_i(t)\big\|\nabla g_i\big(\ox(t)-\ou(t)\big)\big\|\\
&\le\dfrac{\be}{M_1}\Big\|\sum_{i\in I(\ox(t)-\ou(t))}\eta_i(t)\nabla g_i\big(\ox(t)-\ou(t)\big)\Big\|\\
&\le\dfrac{\be}{M_1}\|\dot\ox(t)\|+\dfrac{\be}{M_1}\big\|f\big(\ox(t),\oa(t)\big)\big\|
\end{aligned}
$$
valid for a.e.\ $t\in[0,T]$ and all $i=1,\ldots,m$, which are consequences of \eqref{e:pddr}, \eqref{e:bound-grad}, and \eqref{e:w-inverse-triangle-in}. The uniqueness of $\eta(t)$ a.e.\ on $[0,T]$ follows from the positive linear independence of the gradients $\nabla g_i(x)$ on $C$, which is a consequence of the standing assumptions in Section~2. \\\vspace*{0.02in}
{\bf Step~3:} {\em Constructions of approximating dual elements on $[0,T]$.} The next step is to extend the discrete dual elements from Theorem~\ref{Th:OC-DP} on the continuous-time interval $[0,T]$ in the way appropriate for the subsequent limiting procedure. Define $q^k(t)=(q^{xk}(t),q^{uk}(t),q^{ak}(t))$ on $[0,T]$ as the piecewise linear extensions of $q^k(t^k_j):=p^k_j$ when $j=0,\ldots,k$. Then construct $\gg^k(t)$ on $[0,T]$ by
$$
\gg^k(t):=\gg^k_j\;\mbox{ for }\;t\in[t^k_j,t^k_{j+1}),\;j=0,\ldots,k-1,\;\mbox{ with }\;\gg^k(t^k_k):=0.
$$
and define further $\xi^{1k}(t):=\dfrac{\xi^{1k}_j}{h_k},\;\xi^{2k}(t):=\dfrac{\xi^{2k}_j}{h_k}$ for $t\in[t^k_j,t^k_{j+1})$ and $j=0,\ldots,k-1$ with $\xi^{1k}(t^k_k):=\xi^{2k}_k$ and $\xi^{1k}(t^k_k):=\xi^{2k}_k$. It follows from the relationships in \eqref{e:dac16}--\eqref{e:dac18} with $\nu^k(t)$ given in \eqref{max-time} that
\begin{equation*}
\begin{aligned}
\dot q^{xk}(t)-\lm^kw^{xk}(t)&=\nabla_xf\big(\ox^k(\nu^k(t)),\oa^k(\nu^k(t))\big)^*\big(\lm^k(v^{xk}(t)+\th^{xk}(t)\big)-q^{xk}\big(\nu^k(t)+h_k)\big)\\
&-\sum^m_{i=1}\eta^k_i(t)\nabla^2g_i\big(\ox^k(\nu^k(t))-\ou^k(\nu^k(t))\big)\big(\lm^k(v^{xk}(t)+\th^{xk}(t)\big)-q^{xk}\big(\nu^k(t)+h_k)\big)\\
&-\sum^m_{i=1}\gg^k_i(t)\nabla g_i\big(\ox^k(\nu^k(t))-\ou^k(\nu^k(t))\big),
\end{aligned}
\end{equation*}
\begin{equation*}
\begin{aligned}
\dot q^{uk}(t)-\lm^kw^{uk}(t)&=\sum^m_{i=1}\eta^k_i(t)\nabla^2g_i\big(\ox^k(\nu^k(t))-\ou^k(\nu^k(t))\big)\big(\lm^k(v^{xk}(t)+\th^{xk}(t)\big)-q^{xk}\big(\nu^k(t)+h_k)\big)\\
&+\sum^m_{i=1}\gg^k_i(t)\nabla g_i\big(\ox^k(\nu^k(t))-\ou^k(\nu^k(t))\big)+2\big(\xi^{1k}(t)+\xi^{2k}(t)\big)\ou^k(\nu^k(t)\big),
\end{aligned}
\end{equation*}
\begin{equation*}
\dot q^{ak}(t)-\lm^kw^{ak}(t)=\nabla_af\big(\ox^k(\nu^k(t)),\oa^k(\nu^k(t))\big)^*\big(\lm^k(v^{xk}(t)+\th^{xk}(t)\big)-q^{xk}\big(\nu^k(t)+h_k)\big)
\end{equation*}
for $t\in(t^k_j,t^k_{j+1})$ and $j=0,\ldots,k-1$. The next triple is $p^k(t)=(p^{xk}(t),p^{uk}(t),p^{ak}(t))$ defined by
\begin{equation*}
\begin{aligned}
p^k(t):=q^k(t)+\int_{[t,T]}&\bigg(-\sum^m_{i=1}\eta^k_i(t)\nabla^2g_i\big(\ox^k(\nu^k(t))-\ou^k(\nu^k(t))\big)\big(\lm^k(v^{xk}(t)+\th^{xk}(t)\big)-q^{xk}\big(\nu^k(t)+h_k)\big)\\
&-\sum^m_{i=1}\gg^k_i(s)\nabla g_i\big(\ox^k(\nu^k(s))-\ou^k(\nu^k(s))\big),\\
&\sum^m_{i=1}\eta^k_i(t)\nabla^2g_i\big(\ox^k(\nu^k(t))-\ou^k(\nu^k(t))\big)\big(\lm^k(v^{xk}(t)+\th^{xk}(t))-q^{xk}(\nu^k(t)+h_k))\\
&+\sum^m_{i=1}\gg^k_i(t)\nabla g_i(\ox^k(\nu^k(t))-\ou^k(\nu^k(t)))+2\big(\xi^{1k}(t)+\xi^{2k}(t)\big)\ou^k\big(\nu^k(t)\big),0\bigg)\,ds
\end{aligned}
\end{equation*}
for all $t\in[0,T]$. We clearly have $p^k(T)=q^k(T)$ together with the differential condition
\begin{equation*}
\begin{aligned}
\dot p^k(t)=\dot q^k(t)-&\bigg(-\sum^m_{i=1}\eta^k_i(t)\nabla^2g_i\big(\ox^k(\nu^k(t))-\ou^k(\nu^k(t))\big)\big(\lm^k(v^{xk}(t)+\th^{xk}(t)\big)-q^{xk}\big(\nu^k(t)+h_k)\big)\\
&-\sum^m_{i=1}\gg^k_i(t)\nabla g_i\big(\ox^k(\nu^k(t))-\ou^k(\nu^k(t))\big),\\
&\sum^m_{i=1}\eta^k_i(t)\nabla^2g_i\big(\ox^k(\nu^k(t))-\ou^k(\nu^k(t))\big)\big(\lm^k(v^{xk}(t)+\th^{xk}(t)\big)-q^{xk}\big(\nu^k(t)+h_k)\big)\\
&+\sum^m_{i=1}\gg^k_i(t)\nabla g_i\big(\ox^k(\nu^k(t))-\ou^k(\nu^k(t))\big)+2\big(\xi^{1k}(t)+\xi^{2k}(t)\big)\ou^k\big(\nu^k(t)\big),0\bigg)
\end{aligned}
\end{equation*}
valid for a.e.\ $t\in[0,T]$. It follows from the relationships above that
\begin{equation}\label{e:px-der}
\begin{aligned}
\dot p^{xk}(t)-\lm^kw^{xk}(t)&=\nabla_xf\big(\ox^k(\nu^k(t)),\oa^k(\nu^k(t))\big)^*\big(\lm^k(v^{xk}(t)+\th^{xk}(t)\big)-q^{xk}\big(\nu^k(t)+h_k)\big),
\end{aligned}
\end{equation}
\begin{equation}\label{e:pu-der}
\dot p^{uk}(t)-\lm^kw^{uk}(t)=0,
\end{equation}
\begin{equation}\label{e:pa-der}
\dot p^{ak}(t)-\lm^kw^{ak}(t)=\nabla_af\big(\ox^k(\nu^k(t)),\oa^k(\nu^k(t))\big)^*\big(\lm^k(v^{xk}(t)+\th^{xk}(t)\big)-q^{xk}\big(\nu^k(t)+h_k)\big),
\end{equation}
for every $t\in(t^k_j,t^k_{j+1}),\;j=0,\ldots,k-1$. Now we get the measures $\gg^k$, $\xi^{1k}$, and $\xi^{2k}$ on $[0,T]$ given by
\begin{equation}\label{e:discrete-measures}
\left\{
\begin{array}{ll}
\disp\int_Ad\gg^k:=\int_A\Big(\sum^m_{i=1}\eta^k_i(t)\nabla^2g_i\big(\ox^k(\nu^k(t))-\ou^k(\nu^k(t))\big)\big(\lm^k(v^{xk}(t)+\th^{xk}(t)\big)\\
-q^{xk}\big(\nu^k(t)+h_k)\big)\Big)dt+\disp\int_A\Big(\sum^m_{i=1}\gg^k_i(t)\nabla g_i\big(\ox^k(\nu^k(t))-\ou^k(\nu^k(t))\big)\Big)dt,\\
\disp\int_Ad\xi^{1k}(t)dt:=\int_A\xi^{1k}(t)dt,\;\int_Ad\xi^{2k}(t)dt:=\int_A\xi^{2k}(t)dt
\end{array}\right.
\end{equation}
for any Borel subset $A\subset[0,T]$. Finally in the step, we employ the standard normalization procedure to equivalently rewrite the nontriviality condition \eqref{e:dac26} in the form
\begin{equation}\label{e:discrete-nontri}
\lm^k+\left\|p^k(T)\right\|+\left\|q^{uk}(0)\right\|+\int^T_0\left|\xi^{1k}(t)\right|dt+\int^T_0\left|\xi^{2k}(t)\right|dt+|\xi^{1k}_k|+|\xi^{2k}_k|=1,\;\;k\in\N.
\end{equation}
{\bf Step~4:} {\em Verifying the dual dynamic conditions.} By \eqref{e:discrete-nontri} we get $\lm^k\to\lm\ge0$ along a subsequence. To verity next that $\{q^{uk}(\cdot)\}$ is of uniformly bounded variations on $[0,T]$, observe by \eqref{e:dac5} that
$$
\begin{aligned}
\sum^{k-1}_{j=0}\|q^{uk}(t_{j+1})-q^{uk}(t_j)\|&=\sum^{k-1}_{j=0}\|p^{uk}_{j+1}-p^{uk}_j\|\le\sum^{k-1}_{j=1}\|p^{uk}_{j+1}-p^{uk}_j\|+\|p^{uk}_1\|+\|p^{uk}_0\|\\
&\le\lm^k\sum^{k-1}_{j=1}\left\|\dfrac{\th^{uk}_j-\th^{uk}_{j-1}}{h_k}\right\|+\lm^k\dfrac{\|\th^{uk}_0\|}{h_k}+\lm^k\sum^{k-1}_{j=1}\|v^{uk}_j-v^{uk}_{j-1}\|+\lm^k\|v^{uk}_0\|
+\|p^{uk}_0\|\\
&\le2\lm^k\sum^{k-1}_{j=1}\left\|\dfrac{\ou^k_{j+1}-2\ou^k_j+\ou^k_{j-1}}{h_k}\right\|+2\lm^k\sum^{k-1}_{j=1}\left\|\dfrac{\ou^k(t^k_{j+1})-2\ou^k(t^k_j)+\ou^k(t^k_{j-1})}{h_k}\right\|\\
&+\lm^k\dfrac{\|\th^{uk}_0\|}{h_k}+\lm^k\sum^{k-1}_{j=1}\|v^{uk}_j-v^{uk}_{j-1}\|+\lm^k\|v^{uk}_0\|+\|p^{uk}_0\|\\
&\le 4\Tilde\mu\lm^k+\lm^k\dfrac{\|\th^{uk}_0\|}{h_k}+\lm^k\sum^{k-1}_{j=1}\|v^{uk}_j-v^{uk}_{j-1}\|+\lm^k\|v^{uk}_0\|+\|p^{uk}_0\|
\end{aligned}
$$
with $\Tilde\mu$ taken from Theorem~\ref{Th:DA}. The differentiability of $\ell_2$ in \eqref{separ} with respect to $\dot u$ yields
$$
v^{uk}_j=\nabla_{\dot u}\ell_2\left(t_j,\dfrac{\ou^k_{j+1}-\ou^k_j}{h_k},\dfrac{\oa^k_{j+1}-\oa^k_j}{h_k}\right)\;\mbox{for}\;\;j=0,\ldots,k-1.
$$
Then the third estimate above ensures that
$$
\lm^k\sum^{k-1}_{j=1}\|v^{uk}_j-v^{uk}_{j-1}\|\le \lm^k\sum^{k-1}_{j-1}L\left(t^k_{j+1}-t^k_j+\left\|\dfrac{\ou^k_{j+1}-2\ou^k_j+\ou^k_{j-1}}{h_k}\right\|\right)\le\lm^kL(T+\Tilde\mu),
$$
which in turn yields by the construction of $\th^{uk}$ in \eqref{e:dac9} the relationships
$$
\dfrac{\th^{uk}_0}{h_k}=\dfrac{2(\ou^k_1-\ou^k_0)}{h_k}-\dfrac{2(\ou(h_k)-\ou(0))}{h_k}\;\mbox{ and }\;
\lm^k\dfrac{\|\th^{uk}_0\|}{h_k}\le 2\lm^k(\mu+\Tilde\mu)
$$
due to the first estimate in \eqref{e:discrete-estimate}. Furthermore
$$
\|v^{uk}_0\|=\left\|\nabla_{\dot u}\ell_2\left(0,\dfrac{\ou^k_1-\ou^k_0}{h_k},\dfrac{\oa^k_1-\oa^k_0}{h_k}\right)\right\|\le L\|\dfrac{\ou^k_1-\ou^k_0}{h_k}\|\le L\Tilde\mu
$$
due to the first estimate in \eqref{estimate1}. This tells us that
\begin{equation}
\label{e:qu-bounded-variation}
\begin{aligned}
\sum^{k-1}_{j=0}\|q^{uk}(t_{j+1})-q^{uk}(t_j)\|&=\sum^{k-1}_{j=0}\|p^{uk}_{j+1}-p^{uk}_j\|\\
&\le 4\Tilde\mu\lm^k+2\lm^k(\mu+\Tilde\mu)+\lm^kL(T+\Tilde\mu)+\lm^kL\Tilde\mu+\|p^{uk}_0\|\\
&=\lm^k(2\mu+6\Tilde\mu+LT+L\Tilde\mu)+\|p^{uk}_0\|
\end{aligned}
\end{equation}
and verifies therefore the uniform bounded variations of the sequence $\{q^{uk}(\cdot)\}$ on $[0,T]$.

Our next goal is prove the boundedness of $\{(p^{xk}_0,\ldots,p^{xk}_k)\}$. It follows from \eqref{e:dac16} and \eqref{e:dac17} that
\begin{eqnarray}\label{e:qxk-variation}
\begin{array}{ll}
p^{xk}_{j+1}-p^{xk}_j&=h_k\lm^k(w^{xk}_j+w^{uk}_j)+2(\xi^{1k}_j+\xi^{2k}_j)\ou^k_j\\
&+h_k\nabla_xf(\ox^k_j,\oa^k_j)^*(\lm^k(v^{xk}_j+h^{-1}_k\th^{xk}_j)-p^{xk}_{j+1})-(p^{uk}_{j+1}-p^{uk}_j),
\end{array}
\end{eqnarray}
and readily implies the estimates
\begin{equation}\label{e:pxk-estimate}
\begin{aligned}
\|p^{xk}_j\|&\le(1+h_k\|\nabla_xf(\ox^k_j,\oa^k_j)\|)\|p^{xk}_{j+1}\|+h_k\lm^k(\|w^{xk}_j\|+\|w^{uk}_j\|)+2r_2(\xi^{1k}_j+\xi^{2k}_j)\\
&+\|\nabla_xf(\ox^k_j,\oa^k_j)\|\|h_kv^{xk}_j\|+h_k\|\nabla_xf(\ox^k_j,\oa^k_j)\|\lm^k\|\th^{xk}(t_j)\|+\|p^{uk}_{j+1}-p^{uk}_j\|
\end{aligned}
\end{equation}
valid for all $j=0,\ldots,k-1$. Denoting further
$$
A^k_j:=h_k\lm^k(\|w^{xk}_j\|+\|w^{uk}_j\|)+2r_2(\xi^{1k}_j+\xi^{2k}_j)+\|\nabla_xf(\ox^k_j,\oa^k_j)\|\|h_kv^{xk}_j\|
+h_k\|\nabla_xf(\ox^k_j,\oa^k_j)\|\lm^k\|\th^{xk}(t_j)\|+\|p^{uk}_{j+1}-p^{uk}_j\|
$$
and selecting a $\Tilde{M_1}$ such that $\|\nabla_xf(\ox^k_j,\oa^k_j)\|\le\tilde{M_1}$ for all $j=0,\ldots,k-1$ and $k\in\N$, we have by \eqref{e:discrete-nontri}
$$
\begin{aligned}
\sum^{k-1}_{j=0}h_k\left\|\nabla_xf(\ox^k_j,\oa^k_j)\right\|\lm^k\left\|\th^{xk}(t_j)\right\|
&\le\lm^k\Tilde{M_1}\sum^{k-1}_{j=0}\sqrt{h_k\int^{t_{j+1}}_{t_j}\left\|\th^{xk}(t)\right\|^2dt}\\
&\le \lm^k\Tilde{M_1}\sqrt{\int^T_0\left\|\th^{xk}(t)\right\|^2dt}\dn0\;\mbox{as}\;k\to\infty.
\end{aligned}
$$
The structure \eqref{separ} of the running cost in \eqref{e:Bolza} and the imposed assumptions on its Lipschitz constant $L(t)$  yield by \eqref{subg} the relationships
$$
\sum^{k-1}_{j=0}\left\|h_kw^{xk}_j\right\|=\sum^{k-1}_{j=0}h_k\left\|w^{xk}(t_j)\right\|\le\sum^{k-1}_{j=0}L(t_j)h_k\le 2\int^T_0L(t)dt:=\Tilde L<\infty,
$$
which ensure furthermore that $\disp\sum^{k-1}_{j=0}\left\|h_kw^{uk}_j\right\|\le\Tilde L$ and $\disp\sum^{k-1}_{j=0}\left\|h_kv^{xk}_j\right\|\le\Tilde L$. We also have
$$
\sum^{k-1}_{i=0}|\xi^{1k}_j+\xi^{2k}_j|=\int^T_0|\xi^{1k}(t)+\xi^{2k}(t)|dt\le 1.
$$
As follows from the above arguments, the boundedness of $\{\|\nabla_xf(\ox^k_j,\oa^k_j)\|\}$ with the usage of \eqref{e:qu-bounded-variation} that
\begin{equation}\label{estimate1a}
\sum^{k-1}_{j=0}A^k_j\le \Tilde{M_2}
\end{equation}
for some constant $\Tilde{M_2}>0$. Combining it with \eqref{e:pxk-estimate} gives us the estimates
\begin{equation*}
\|p^{xk}_j\|\le(1+\Tilde{M_1}h_k)\|p^{xk}_{j+1}\|+A^k_j\;\mbox{ for all }\;j=0,\ldots,k-1.
\end{equation*}
Proceeding now by induction shows that
$$
\begin{aligned}
\|p^{xk}_j\|&\le(1+\Tilde{M_1}h_k)^{k-j}\|p^{xk}_k\|+\sum^{k-1}_{i=j}A^k_i(1+\Tilde{M_1}h_k)^{i-j}\\
&\le e^{\Tilde{M_1}}+e^{\Tilde{M_1}}\sum^{k-1}_{i=0}A^k_i\le e^{\Tilde{M_1}}(1+\Tilde{M_2})
\end{aligned}
$$
for $j=0,\ldots,k-1$, which justifies the boundedness of $\{(p^{xk}_0,\ldots,p^{xk}_k)\}$.\vspace*{-0.02in}

Next we show that the functional sequences $\{q^{xk}(\cdot)\}$ and $\{q^{ak}(\cdot)\}$ are of uniform bounded variations on $[0,T]$. It follows from \eqref{e:qxk-variation} that
$$
\sum^{k-1}_{j=0}\left\|q^{xk}(t_{j+1})-q^{xk}(t_j)\right\|=\sum^{k-1}_{j=0}\left\|p^{xk}_{j+1}-p^{xk}_j\right\|\le\sum^{k-1}_{j=0}A^k_j
+\sum^{k-1}_{j=0}h_k\left\|\nabla_xf(\ox^k_j,\oa^k_j)\right\|\left\|p^{xk}_{j+1}\right\|,
$$
which verifies the claimed property for $\{q^{xk}(\cdot)\}$ due to \eqref{estimate1a} and the boundedness of $\{(p^{xk}_0,\ldots,p^{xk}_k)\}$. Furthermore, \eqref{e:dac18} leads us to the estimate
$$
\sum^{k-1}_{j=0}\left\|q^{ak}(t_{j+1})-q^{ak}(t_j)\right\|\le h_k\sum^{k-1}_{j=0}\left\|\lm^kw^{ak}_j\right\|+h_k\sum^{k-1}_{j=0}\left\|\nabla_af(\ox^k_j,\oa^k_j)\right\|\left(\lm^k\left(\left\|v^{xk}_j\right\|+\left\|\th^{xk}(t^k_j)\right\|
\right)+\left\|p^{xk}_{j+1}\right\|\right),
$$
which ensures the uniform bounded variations of $\{q^{ak}(\cdot)\}$ and hence of the whole sequence of triples $\{q^k(\cdot)\}$. We clearly get the validity of
$$
2\left\|q^k(t)\right\|-\left\|q^k(0)\right\|-\left\|q^k(T)\right\|\le\left\|q^k(t)-q^k(0)\right\|+\left\|q^k(T)-q^k(t)\right\|\le\var\left(q^k;[0,T]\right)
$$
whenever $t\in[0,T]$. Hence $\{q^k(\cdot)\}$ is bounded on $[0,T]$ due to the boundedness of $\{q^k(0)\}$ and $\{q^k(T)\}$. The classical Helly theorem provides a function of bounded variation $q(\cdot)$ such that $q^k(t)\to q(t)$ for all $t\in[0,T]$. It follows from \eqref{e:discrete-nontri} that the sequences $\{\xi^{1k}\}$ and $\{\xi^{2k}\}$ are bounded in $C^*([0,T];\R_+)$ and $C^*([0,T];\R_-)$, respectively. It allows us to get the boundedness of $\{\gg^k\}$ in $C^*([0,T];\R^n)$ from \eqref{e:dac17} and the uniform bounded variations of $\{q^{uk}(\cdot)\}$ on $[0,T]$. We derive from the weak$^*$ sequential compactness of balls in these spaces that there are $\gg\in C^*([0,T];\R^n)$, $\xi^1\in C^*([0,T];\R_+)$, and $\xi^2\in C^*([0,T];\R_-)$ such that the triples $(\gg^k,\xi^{1k},\xi^{2k})$ weak$^*$ converge to $(\gg,\xi^1,\xi^2)$ along some subsequence.\vspace*{-0.05in}

Using \eqref{e:px-der}--\eqref{e:pa-der} and \eqref{e:discrete-nontri} together with the uniform boundedness of $q^k(\cdot),\;w^k(\cdot)$, and $v^k(\cdot)$ on $[0,T]$ ensures the boundedness of the sequence $\{p^k(\cdot)\}$ in $W^{1,2}([0,T];\R^{3n})$ and hence its weak compact in this space. By Mazur's theorem we find a function $p(\cdot)\in W^{1,2}([0,T];\R^{3n})$ such that a sequence of convex combinations of $\dot p^k(t)$ converges to $\dot p(t)$ for a.e.\ $t\in[0,T]$. Then the passage to the limit in \eqref{e:px-der}--\eqref{e:pa-der} justifies the claimed representation of $\dot p(\cdot)$ in \eqref{e:EL}.\vspace*{-0.05in}

Our next aim is to derive the optimality conditions of the theorem that involve the dual arc $q(\cdot)$ of bounded variation on $[0,T]$. Observe that the condition $\eta_i(t)>0$ for some $t\in[0,T]$ and $i\in\{1,\ldots,m\}$ yields $\eta^k_i(t)>0$ for large $k$ by the a.e.\ convergence $\eta^k_i(\cdot)\to\eta_i(\cdot)$ on $[0,T]$. This implies by \eqref{e:dac25} that
$$
\big\la\nabla g_i\big(x^k(t)-u^k(t)\big),-q^{xk}\big(\nu(t)+h_k\big)+\lm^k\big(\th^{xk}(t)+v^{xk}(t)\big)\big\ra=0
$$
for such $k$ and $t$, and hence we get by passing to the limit that
$$
\big\la\nabla g_i\big(\ox(t)-\ou(t)\big),\lm v^x(t)-q^x(t)\big\ra=0,
$$
which verifies the second implication in \eqref{e:implications}. By the construction of $q^k(\cdot)$ we get
\begin{equation}\label{e:quk-qak}
q^{uk}\big(\nu(t)+h_k\big)=\lm^k\big(v^{uk}(t)+\th^{uk}(t)\big)\;\;\mbox{and}\;\;q^{ak}\big(\nu(t)+h_k\big)=\lm^k\big(v^{ak}(t)+\th^{ak}(t)\big)
\end{equation}
whenever $t\in(t^k_j,t^k_{j+1})$ and $j=0,\ldots,k-1$. Involving \eqref{subg} and the assumptions on $\ell_2,\ell_3$ in \eqref{separ} gives us both conditions in \eqref{e:qu-qa} by passing to the limit in \eqref{e:quk-qak}. Proceeding similarly to \cite[p. 325]{v} yields
$$
\begin{aligned}
\bigg\|\int_{[t,T]}&\bigg(\sum^m_{i=1}\eta^k_i(s)\nabla^2 g_i\big(\ox^k(\nu^k(s))-\ou^k(\nu^k(s))\big)\big(\lm^k(v^{xk}(s)+\th^{xk}(s)\big)-q^{xk}\big(\nu^k(s)+h_k)\big)\\
&+\sum^m_{i=1}\gg^k_i(s)\nabla g_i\big(\ox^k(\nu^k(s))-\ou^k(\nu^k(s))\big)\bigg)ds-\int_{[t,T]}d\gg(s)\bigg\|\\
&=\bigg\|\int_{[t,T]}d\gg^k(s)-\int_{[t,T]}d\gg(s)\bigg\|\to 0\;\;\mbox{as}\;\;k\to\infty
\end{aligned}
$$
for all $t\in[0,T]$ except a countable subset of $[0,T]$. It tells us by using \eqref{e:discrete-measures} that
\begin{equation}
\label{e:gg-convergence}\begin{aligned}
\int_{[t,T]}&\bigg(\sum^m_{i=1}\eta^k_i(s)\nabla^2 g_i\big(\ox^k(\nu^k(s))-\ou^k(\nu^k(s))\big)\big(\lm^k(v^{xk}(s)+\th^{xk}(s)\big)-q^{xk}\big(\nu^k(s)+h_k)\big)\\
&+\sum^m_{i=1}\gg^k_i(s)\nabla g_i\big(\ox^k(\nu^k(s))-\ou^k(\nu^k(s))\big)\bigg)ds\to\int_{[t,T]}d\gg(s)\;\;\mbox{as}\;\;k\to\infty.
\end{aligned}
\end{equation}
To derive \eqref{e:p-q relation} by passing to the limit in the differential condition for $p^k(t)$, consider the estimate
\begin{equation}\label{estimate3}
\begin{aligned}
&\left\|\int_{[t,T]}\big(\xi^{1k}(s)+\xi^{2k}(s)\big)\ou^k\big(\nu^k(s)\big)ds-\int_{[t,T]}\ou(s)d(\xi^1(s)+\xi^2(s))\right\|\\
\le&\left\|\int_{[t,T]}\big(\xi^{1k}(s)+\xi^{2k}(s)\big)\ou^k\big(\nu^k(s)\big)ds-\int_{[t,T]}\big(\xi^{1k}(s)+\xi^{2k}(s)\big)\ou(s)ds\right\|\\
+&\left\|\int_{[t,T]}\big(\xi^{1k}(s)+\xi^{2k}(s)\big)\ou(s)ds-\int_{[t,T]}\ou(s)d\big(\xi^1(s)+\xi^2(s)\big)\right\|\\
=&\left\|\int_{[t,T]}(\xi^{1k}(s)+\xi^{2k}(s))\big[\ou^k\big(\nu^k(s)\big)-\ou(s)\big]ds\right\|\\
+&\left\|\int_{[t,T]}\big(\xi^{1k}(s)+\xi^{2k}(s)\big)\ou(s)ds-\int_{[t,T]}\ou(s)d\big(\xi^1(s)+\xi^2(s)\big)\right\|.
\end{aligned}
\end{equation}
Note that the first summand after the equality sign in \eqref{estimate3} vanishes as $k\to\infty$ due to the uniform convergence $\ou^k(\cdot)\to\ou(\cdot)$ on $[0,T]$ and the uniform boundedness of $\disp\int^T_0|\xi^{1k}(t)|dt+\int^T_0|\xi^{2k}(t)|dt$ by \eqref{e:discrete-nontri}. The second summand there also converges to zero for all $t\in[0,T]$ except some countable subset by the weak$^*$ convergence of $\xi^{1k}\to\xi^1$ in $C^*([0,T];\R_+)$ and $\xi^{2k}\to\xi^2$ in $C^*([0,T];\R_-)$. Thus we get
$$
\int_{[t,T]}(\xi^{1k}(s)+\xi^{2k}(s))\ou^k(\tau^k(s))ds\to\int_{[t,T]}\ou(s)d(\xi^1(s)+\xi^2(s))\;\mbox{ as }\;k\to\infty
$$
and hence obtain \eqref{e:p-q relation} by passing to the limit in the above differential condition for $p^k(t)$.\\
{\bf Step~5:} {\em Verifying transversality.} Observe first that the sequence $\{\eta^k_{ki}\}$ admits a convergent subsequence. To show it, we deduce from the second discrete transversality condition in \eqref{e:dac19} that
$$
\sum_{i\in I(\ox(T)-\ou(T))}\eta^k_{ki}\nabla g_i(\ox^k_k-\ou^k_k)=\sum_{i\in I(\ox^k_k-\ou^k_k)}\eta^k_{ki}\nabla g_i(\ox^k_k-\ou^k_k)=-p^{uk}_k-2(\xi^{1k}_k+\xi^{2k}_k)\ou^k_k
$$
with $\eta^k_{ki}=0$ for $i\in\{1,\ldots,m\}\backslash I(\ox^k_k-\ou^k_k)$. This justifies by \eqref{e:discrete-nontri} the boundedness of the sequence
$\disp\left\{\sum_{i\in I(\ox(T)-\ou(T))}\eta^k_{ki}\nabla g_i(\ox^k_k-\ou^k_k)\right\}$. It follows from the standing assumptions in \eqref{e:bound-grad} and \eqref{e:w-inverse-triangle-in} that
$$
\eta^k_{ki}\le\dfrac{1}{M_1}\sum_{j\in I(\ox(T)-\ou(T))}\eta^k_{kj}\left\|\nabla g_i(\ox^k_k-\ou^k_k)\right\|\le\dfrac{\be}{M_1}\left\|\sum_{i\in I(\ox(T)-\ou(T))}\eta^k_{ki}\nabla g_i(\ox^k_k-\ou^k_k)\right\|,
$$
which ensures the boundedness of $\{\eta^k_{ki}\}$ for $i=1,\ldots,m$, and thus we get $\eta^k_{ki}\to\Tilde\eta_i$ along a subsequence of $k\to\infty$. Denote $\eta(T):=(\Tilde\eta_1,\ldots,\Tilde\eta_n)$ and observe that it is well defined due to the positive linear independence of vectors $\nabla g_i(\ox(T)-\ou(T))$, which follows from the standing assumptions. Then
$$
\vTh_k:=\sum_{i\in I(\ox(T)-\ou(T))}\eta^k_{ki}\nabla g_i(\ox^k_k-\ou^k_k)\to\vTh:=\sum_{i\in I(\ox(T)-\ou(T))}\eta_i(T)\nabla g_i\big(\ox(T)-\ou(T)\big)
$$
as $k\to\infty$, where the defined vector $\vTh$ satisfies the inclusion in \eqref{e:right-trans1}. Further, it follows from \eqref{e:dac22a} and the second equation in \eqref{e:dac19} that
\begin{eqnarray}\label{trans1}
p^{uk}_k+\vTh_k=-2(\xi^{1k}_k+\xi^{2k}_k)\in-2\ou^k_k\left(N_{[0,r_2+\ve_k]}(\|\ou^k_k\|)+N_{[r_1-\ve_k,\infty)}(\|\ou^k_k\|)\right).
\end{eqnarray}
Passing now to the limit in \eqref{trans1} and \eqref{e:dac19} with taking into account the subdifferential robustness as well as the convergence of $\{\xi^{1k}_k\}$ and $\{\xi^{2k}_k\}$ \eqref{e:dac22a}, we justify the transversality conditions \eqref{e:right-trans}.\\
{\bf Step~6:} {\em Verifying measure nonatomicity.} We provide the verification of the nonatomicity condition (a) while observing that the case of (b) is similar. To proceed, pick $t\in[0,T]$ with $g_i(\ox(t)-\ou(t))>0$ whenever $i=1,\ldots,m$ and employing the continuity of $g_i$ and $(\ox(\cdot),\ou(\cdot))$ get a neighborhood $V_t$ of $t$ such that $g_i(\ox(s)-\ou(s))>0$ for all $s\in V_t$ and $i=1,\ldots,m$. The obtained convergence of the discrete optimal solutions tells us that $g_i(\ox^k(t^k_j)-\ou^k(t^k_j))>0$ if $t^k_j\in V_t$ for $i=1,\ldots,m$ and all large $k$. It follows from \eqref{e:dac20} and \eqref{e:dac23} that $\eta^k_{ji}=0$ and $\gg^k_{ji}=0$ for $i=1,\ldots,m$. Thus
$$
\|\gg^k\|(V)=\disp\int_Vd\left\|\gg^k\right\|=\int_V\left\|\gg^k(t)\right\|dt=0
$$
by \eqref{e:discrete-measures}. Passing to the limit as $k\to\infty$ with taking into account the measure convergence established in Step~3, we get $\|\gg\|(V)=0$, which verifies the claimed measure nonatomicity.\\
{\bf Step~7:} {\em Verifying nontriviality.} To start with the verification of the general nontriviality condition \eqref{e:nontri-condition}, suppose on the contrary that $\lm=0,\; q^u(0)=0,\;p(T)=0,\;\|\xi^1\|_{TV}=0$, and $\|\xi^2\|_{TV}=0$. Then
$$
\lm^k\to 0,\;q^{uk}(0)\to 0,\;p^k(T)\to 0,\;\disp\int_{[0,T]}|\xi^{1k}(t)|dt\to 0,\;\disp\int_{[0,T]}|\xi^{2k}(t)|dt\to 0
$$
as $k\to\infty$. Let us check that in this case we get $\xi^{1k}_k+\xi^{2k}_k\to 0$. Indeed, observe that the convergence $p^k(T)\to 0,\;\lm^k\to0$ implies by the first condition in \eqref{e:dac19} that $p^{xk}_k\to 0$, $p^{uk}_k\to 0$, and $\disp\sum^m_{i=1}\eta^k_{ki}\nabla g_i(\ox^k_k-\ou^k_k)\to0$. Then the second condition in \eqref{e:dac19} yields $(\xi^{1k}_k+\xi^{2k}_k)\ou^k_k\to 0$, and $\xi^{1k}_k+\xi^{2k}_k\to 0$ due to $\ou^k_k\not=0$ for all large $k$. This clearly contradicts \eqref{e:discrete-nontri}, and so we are done with \eqref{e:nontri-condition}.\vspace*{-0.02in}

It remains to verify the enhanced nontriviality conditions of the theorem under the additional assumptions made. To proceed with \eqref{e:enhanced-nontri}, suppose that $g_i(x_0-\ou(0))>0$ for $i=1,\ldots,m$ and, arguing by contradiction, that $\lm=0,\;p(T)=0,\;\|\xi^1\|_{TV}=0$, and $\|\xi^2\|_{TV}=0$. Then $g_i(\ox^k_0-\ou^k_0)>0$ for large $k$. It follows from \eqref{e:dac20} and \eqref{e:dac21} that $\eta^k_{ji}=0$ and $\gg^k_{ji}=0$ for $i=1,\ldots,m$. Unifying it with \eqref{e:dac17} and the construction of $q^{uk}(\cdot)$ in Step~3 shows that
$$
q^{uk}(0)=p^{uk}_0=p^{uk}_1+2(\xi^{1k}_0+\xi^{2k}_0)\ou^k_0+h_k\lm^kw^{uk}_0
$$
whenever $k\in\N$ is sufficiently large. This yields the estimates
$$
\begin{aligned}
\|q^{uk}(0)\|&\le\lm^k\|v^{uk}_0\|+\lm^k\dfrac{\left\|\th^{uk}_0\right\|}{h_k}+\lm^kh_k\left\|w^{uk}_0\right\|+2r_2\left|\xi^{1k}_0+\xi^{2k}_0\right|\\
&\le\lm^kL\Tilde\mu+2\lm^k(\mu+\Tilde\mu)+\lm^k\Tilde\mu+\int_{[0,T]}|\xi^{1k}(t)|dt+\int_{[0,T]}|\xi^{2k}(t)|dt,
\end{aligned}
$$
which imply in turn that $q^u(0)=\disp\lim_{k\to\infty}q^{uk}(0)=0$ while contradicting the nontriviality condition \eqref{e:nontri-condition}. The verification of the other enhanced nontriviality condition \eqref{e:enhanced-nontri1} is similar.$\h$\vspace*{-0.2in}

\section{Numerical Examples}
\setcounter{equation}{0}\vspace*{-0.1in}

In this section we present two examples, which are related to real-life models while illustrating some special features and applications of the obtained necessary optimality conditions for the controlled sweeping process under consideration. In both examples, the obtained optimality conditions allow us to determine optimal solutions and explicitly calculate their parameters.\vspace*{-0.05in}

The first example addresses a one-dimensional sweeping control model of type $(P)$.\vspace*{-0.13in}

\begin{example}{\bf (optimal control of car motion).}\label{ex1} Consider a car moving towards the traffic light with the initial speed $s=9$m/s (about $20$mi/h). When the car is 250 meters away from the traffic light, the light turns into green color and it lasts for 30 seconds. We need to control the motion of the car in such a way that after 20 seconds it must be as close to the traffic light as possible and the energy used to adjust the speed must be minimized as well; see Figure~1.
\begin{center}
\begin{tikzpicture}

\begin{scope}[scale=0.6]
    \shade[top color=red, bottom color=white, shading angle={135}]
    [draw=black,fill=red!20,rounded corners=1.2ex,very thick] (1.5,.5) -- ++(0,1.3) -- ++(1,0) --  ++(3,0) -- ++(1,-0.3) -- ++(0,-1) -- (1.5,.5) -- cycle;
   \draw[very thick, rounded corners=0.5ex,fill=black!20!blue!20!white,thick]  (2.5,1.8) -- ++(0.6,0.7) -- ++(1.4,0) -- ++(0.9,-0.7) -- (2.5,1.8);
    \draw[thick]  (3.8,1.8) -- (3.8,2.5);
    \draw[draw=black,fill=gray!50,thick] (2.75,.5) circle (.5);
    \draw[draw=black,fill=gray!50,thick] (5.5,.5) circle (.5);
    \draw[draw=black,fill=gray!80,semithick] (2.75,.5) circle (.4);
    \draw[draw=black,fill=gray!80,semithick] (5.5,.5) circle (.4);
  \end{scope}
\draw node[below] at (2,0) {-250};
\draw [thick, red](0,0)--(10,0);
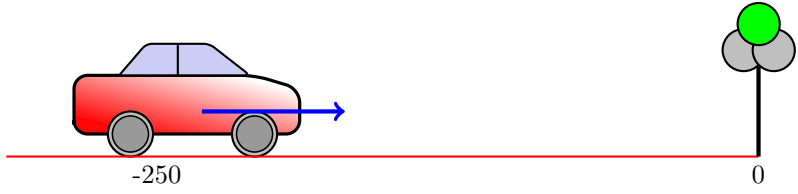
\draw node[below] at (5,-0.5) {{\bf Figure 1:} Direction of optimal control};
\draw [ultra thick, blue] [->] (2.6,0.6)--(4.5,0.6);
\draw node[above] at (10,1) {\TlGreen[4]};
\draw node[below] at (10,0) {0};
\draw [ultra thick] (10,0)--(10,1.2);
\end{tikzpicture}
\end{center}

To solve this problem numerically, let us specify the initial data in problem $(P)$ as follows:
\begin{equation}\label{e:data1}
\left\{\begin{array}{ll}
n=m=d=1,\;T=20,\; x_0:=-250,\;g_1(x):=-x,\;f(x,a):=sa=9a,\\
r_1:=10^{-3},\;r_2:=50,\;\ph:=\dfrac{x^2}{2},\;\mbox{ and }\;\ell(t,x,u,a,\dot x,\dot u,\dot a):=\frac{1}{2}a^2.
\end{array}\right.
\end{equation}
For definiteness, suppose that the traffic light is located at the origin. It follows from \eqref{e:data1} that the set of $a$-controls can be assumed to be uniformly bounded. Thus $(P)$ admits an optimal solution $(\ox(\cdot),\ou(\cdot),\oa(\cdot))\in W^{1,2}([0,20];\R^3)$ by Theorem~\ref{Th:existence}. It is also easy to see that all the assumptions of Theorem~\ref{Th:NOC} are satisfied. Supposing further that $\ox(t)\in\mbox{int}(C+\ou(t))$ for any $t\in[0,20)$ and that $\ox(20)-\ou(20)\in\bd C$, we see that these assumptions are realized for the optimal solution found via the necessary optimality conditions of Theorem~\ref{Th:NOC}. Since we expect the car to be as close to the traffic light as possible after 20 seconds with $\ox(20)=\ou(20)$, the value $\ou(20)$ should be small. The choice of $r_1$ and $r_2$ in \eqref{e:data1} ensures that the validity of the constraints $r_1\le|\ou(t)|\le r_2$ for all $t\in[0,20]$.\vspace*{-0.02in}

Applying the necessary optimality conditions of Theorem~\ref{Th:NOC} gives us the following relationships with a number $\lm\ge0$ and a function $\eta(\cdot)\in L^2([0,20];\R_+)$ well defined at $t=20$:\vspace*{-0.1in}
\begin{enumerate}
\item $w(t)=(0,0,\oa(t)),\;v(t)=(0,0,0)$ a.e.\ $t\in[0,20]$;\vspace*{-0.1in}
\item $\ox(t)<\ou(t)\Longrightarrow \eta(t)=0$ a.e.\ $t\in[0,20]$;\vspace*{-0.1in}
\item $\eta(t)>0\Longrightarrow q^x(t)=0$ a.e.\ $t\in[0,20]$ including $t=20$;\vspace*{-0.1in}
\item $\dot\ox(t)+9\oa(t)=-\eta(t)$ a.e.\ $t\in[0,20]$;\vspace*{-0.1in}
\item $(\dot p^x(t),\dot p^u(t),\;\dot p^a(t))=(0,0,\lm\oa(t)-9q^x(t))$ a.e.\ $t\in[0,20]$;\vspace*{-0.1in}
\item $q^u(t)=0,\;q^a(t)=0$ a.e.\ $t\in[0,20]$;\vspace*{-0.1in}
\item $\disp(q^x(t),q^u(t),q^a(t))=(p^x(t),p^u(t),p^a(t))-\int_{[t,20]}(-d\gg,2\ou(s)d(\xi^1(s)+\xi^2(s)),0)$ a.e.\ $t\in[0,20]$;\vspace*{-0.15in}
\item $-p^x(20)-\eta(20)=\lm\ox(20)$;\vspace*{-0.1in}
\item $p^u(20)+\eta(20)\in-2\ou(20)\left(N_{[0,1]}(|\ou(20)|)+N_{[10^{-3},\infty)}(|\ou(20)|)\right)$;\vspace*{-0.1in}
\item $-\eta(20)\in N_C(\ox(20)-\ou(20))$;\vspace*{-0.1in}
\item $\lm+|q^u(0)|+|p(20)|+\|\xi^1\|_{TV}+\|\xi^2\|_{TV}>0$.\vspace*{-0.1in}
\end{enumerate}\vspace*{-0.02in}
We get from (5)--(7) that $p^x(\cdot)$ is a constant function on $[0,20]$ and that
$$
\lm\oa(t)=9q^x(t)=9p^x(20)+9\int_{[t,20]}d\gg.
$$
Proceeding similarly to \cite[Example~1]{cm2} gives us the relationship
\begin{equation}\label{e:data2}
\lm\oa(t)=9p^x(20)+9\gg(\{20\})\;\;\mbox{ a.e. }\;\;[0,20].
\end{equation}
If $\lm>0$, we can deduce from \eqref{e:data2} that $\oa(\cdot)$ must be a constant function, $\oa(\cdot)\equiv\vTh$, on $[0,20]$ due to its continuity on this interval. In the case $\lm=0$, we may assume that we control the speed of the car at the initial time as $\dot\ox(0)=-10\vTh$ and maintain this speed till the end of the process, i.e., $\dot\ox(t)=-10\vTh$ for a.e.\ $t\in[0,20]$. This results in $\oa(\cdot)\equiv\vTh$ on $[0,20]$. Then (2) and (4) yield
$$
\ox(t)=x_0+\int^t_0\dot\ox(s)ds=-250-9\vTh t\;\;\mbox{for all}\;\;t\in[0,20].
$$
Consequently, the cost functional in our problem $(P)$ is computed as
$$
J[\ox,\ou,\oa]=\dfrac{(-250-180\vTh)^2}{2}+\dfrac{20\vTh^2}{2}
$$
and clearly achieves its absolute minimum at $\vTh=-\frac{45000}{32420}\approx-1.388$. Thus in this case we arrive by the necessary optimality conditions of Theorem~\ref{Th:NOC} at the (local) optimal solution written as
$$
\ox(t)=-250+12.492t,\;\;\oa(t)=-1.388\;\;\mbox{on}\;\;[0,20]
$$
with  $\ou(\cdot)$ being an absolutely continuous function on $[0,20]$ such that $\ou(20)=\ox(20)=-0.16$. This tells us that at the moment when the car is 250 meters away from the traffic light its speed should be switched to $1.388\times9=12.492$m/s, and we should maintain this speed till the end of the process. After 20 seconds, the car is just 0.16 meter away from the traffic light. As seen, the value $12.492$m/s is very close to $12.5=\frac{250}{20}$. To save more energy needed to get such speed, it makes sense to adjust the running cost as follows: $\ell(t,x,u,a,\dot x,\dot u,\dot a):=\frac{100}{2}a^2$. In this case the cost functional is represented by
$$
J[\ox,\ou,\oa]=\dfrac{(-250-180\vTh)^2}{2}+\dfrac{2000\vTh^2}{2}
$$
and achieves its absolute minimum at $\vTh=-\frac{45000}{34400}\approx-1.308$ by using the optimal solution
$$
\ox(t)=-250+11.772t,\;\oa(t)=-1.308,
$$
where the optimal control $\ou(\cdot)$ is absolutely continuous on $[0,20]$ and such that $\ou(20)=\ox(20)=-14.56$. Hence we should switch the speed of the car to $11.772$m/s when it is 250 meters away from the traffic light and maintain this new constant speed till the end of the process. In this way the car is about 14.56 meters away from the traffic light after 20 seconds.\vspace*{-0.02in}
\end{example}

The next example concerns a rather particular case of the two-dimensional crowd motion model that is formalized as as a nonconvex controlled sweeping process of type $(P)$. A more general controlled crowd motion model on the plane is the subject of our adjacent paper \cite{cm-crowd}. Note that the given example demonstrates the usefulness for calculating optimal solutions of the necessary optimality conditions from Theorem~\ref{Th:NOC} with the general nontriviality condition \eqref{e:nontri-condition}.\vspace*{-0.13in}
\begin{example}{\bf (case for optimal control of the planar crowd motion model).}\label{ex:2} We refer the reader to \cite{mv2,ve} for describing an uncontrolled microscopic version of the crowd motion model as a sweeping process. Here we introduce (following the previous corridor version in \cite{cm2}) controls entering both the moving set and perturbations.
In what follows we restrict ourselves to the case of two participants identified with rigid disks of the same radius $R=3$. The center of the $i$-th disk is denoted by $x_i\in\R^2$. To fulfill the nonoverlapping condition crucial in the crowd motion model, the vector of positions $x=(x_1,x_2)\in\R^4$ has to belong to the nonconvex set of feasible configurations (see Figure~2) defined by
$$
C:=\left\{x=(x_1,x_2)\in\R^4\big|\;g(x):=\|x_1-x_2\|-2R\ge 0\right\}.
$$
Suppose that the initial positions of the two participants are
$$
x_1(0):=\left(-48-\dfrac{6}{\sqrt2},48+\dfrac{6}{\sqrt2}\right),\;\;x_2(0):=(-48,48),
$$
and that the exit is located at the origin $(0,0)$.
\begin{center}
\begin{tikzpicture}
\draw [thick, red](-5,0)--(5,0);
\draw [thick, red] (0,0)--(0,6);
\draw [ultra thick, fill = orange] (-0.5,1) rectangle (0.5,0);
\draw node[below] at (0,0){Exit};
\draw node[below] at (0,-1){{\bf Figure 2}};
\draw [thick] (-4,4) circle (0.5);
\draw node[below] at (-4,4){$x_1$};
\draw [thick] (-3.29,3.29) circle (0.5);
\draw node[above] at (-3.29,3.29){$x_2$};
\draw [ultra thick, blue] [->] (-4,4) -- (-3.2,3.2);
\draw [ultra thick, blue] [->] (-3.29,3.29) -- (-2.49,2.49);
\coordinate[label=above left:] (A) at (-4,4);
\coordinate[label=above left:] (B) at (-5,5);
\coordinate[label=above left:] (C) at (-2.5,4);
\draw[->,purple] (A)--(B);
\draw[->,purple] (A)--(C);
\markangle{B}{A}{C}{$135^\circ$}{2}
\end{tikzpicture}
\end{center}
Assume also that the participants exhibit the same behavior and aim to reach the exit by the shortest path. To regulate the participant speeds under the nonoverlapping condition, we use control functions in the moving set
\begin{eqnarray}\label{u-cont}
u_1(t)=u_2(t)\;\;\mbox{for all}\;\;t\in[0,6]\;\mbox{ with the bounds }\;r_1=1\;\mbox{ and }\;r_2=10
\end{eqnarray}
as well as in the perturbations $a(\cdot)=(a_1(\cdot),a_2(\cdot)):[0,20]\to\R^2$ entering via the velocity function
$$
f(x,a):=\left(\dfrac{s_1a_1}{\|x_1\|}x_1,\dfrac{s_2a_2}{\|x_2\|}x_2\right)=\big(s_1a_1\cos\th_1,s_1a_1\sin\th_1,s_2a_2\cos\th_2,s_2a_2\sin\th_2\big),
$$
where $s_1=6$, $s_2=3$, and $\th_1=\th_2=135^\circ$. Then the controlled dynamics is:
$$
\left\{\begin{array}{ll}
-\dot x(t)\in N_{C(t)}\big(x(t))+f(x(t),a(t)\big)\;\;\mbox{a.e.}\;\;t\in[0,6],\\
C(t):=C+u(t),\;u_1(t)=u_2(t),\;\;r_1\le\|u(t)\|\le r_2\;\mbox{on}\;\;$[0,6]$,\;x(0)=x_0\in C(0).
\end{array}\right.
$$
We consider the cost functional given of the form
$$
J[x,u,a]:=\dfrac{1}{2}\left(\|x(6)\|^2+\int^6_0\|a(t)\|^2dt\right)
$$
the meaning of which is to minimize the distance of the two participants to the exit together with the energy of feasible controls $a(\cdot)$ after six seconds. Applying now the necessary conditions of Theorem~\ref{Th:NOC} yields the following, where $\lm\ge0$ and $\eta_{12}\in L^2([0,6];\R_+)$ are well defined at $t=6$:\vspace*{-0.1in}
\begin{enumerate}
\item $w(t)=\big(0,0,\oa(t)\big),\;v(t)=(0,0,0)$ a.e.\ $t\in[0,6]$;\vspace*{-0.1in}

\item $\dot\ox(t)+f\big(\ox(t),\oa(t)\big)=\eta_{12}(t)\nabla
    g\big(\ox(t)-\ou(t)\big)=\left(-\eta_{12}(t)\dfrac{\ox_2(t)-\ox_1(t)}{\|\ox_2(t)-\ox_1(t)\|},\eta_{12}(t)\dfrac{\ox_2(t)-\ox_1(t)}{\|\ox_2(t)-\ox_1(t)\|}\right)$;\vspace*{-0.1in}

\item $\|\ox_2(t)-\ox_1(t)\|>2R\Longrightarrow\eta_{12}(t)=0$ a.e.\ $t\in[0,6]$;\vspace*{-0.1in}

\item $\eta_{12}(t)>0\Longrightarrow\la q^x_2(t)-q^x_1(t),\ox_2(t)-\ox_1(t)\ra=0$ a.e.\ $t\in[0,6]$;\vspace*{-0.1in}

\item $\dot
p(t)=\left(0,0,\lm\oa_1(t)-6\left(-\frac{\sqrt2}{2}q^x_{11}(t)+\frac{\sqrt2}{2}q^x_{12}(t)\right),\lm\oa_2(t)-3\left(-\frac{\sqrt2}{2}q^x_{21}(t)+\frac{\sqrt2}{2}q^x_{22}(t)\right)
\right)$ a.e.\ $t\in[0,6]$;\vspace*{-0.1in}

\item $q^x(t)=p^x(t)+\gg([t,6])$ a.e.\ $t\in[0,6]$;\vspace*{-0.1in}

\item $q^u(t)=p^u(t)-\disp\int_{[t,6]}2\ou(s)d\left(\xi^1(s)+\xi^2(s)\right)+d\gg(s)=0$ a.e.\ $t\in[0,6]$;\vspace*{-0.1in}

\item $q^a(t)=p^a(t)=0$ a.e.\ $t\in[0,6]$;\vspace*{-0.1in}

\item $p^x(6)+\lm\ox(6)=\left(-\eta_{12}(6)\dfrac{\ox_2(6)-\ox_1(6)}{\|\ox_2(6)-\ox_1(6)\|},\eta_{12}(6)\dfrac{\ox_2(6)-\ox_1(6)}{\|\ox_2(6)-\ox_1(6)\|}\right)$;\vspace*{-0.1in}

\item $\begin{array}{ll}
p^u(6)-\left(-\eta_{12}(6)\dfrac{\ox_2(6)-\ox_1(6)}{\|\ox_2(6)-\ox_1(6)\|},\eta_{12}(6)\dfrac{\ox_2(6)-\ox_1(6)}{\|\ox_2(6)-\ox_1(6)\|}\right)
\in-2\ou(6)\left(N_{[0,2]}(\|\ou(6)\|)+N_{[1,\infty]}(\|\ou(6)\|)\right);
\end{array}$\vspace*{-0.1in}
\item $\lm+\|q^u(0)\|+\|p(6)\|+\|\xi^1\|_{TV}+\|\xi^2\|_{TV}>0$.
\end{enumerate}\vspace*{-0.1in}
Since the two participants are in contact at the initial time, i.e., $\|\ox_2(0)-\ox_1(0)\|=6$, they have the same velocity and maintain it till the end of process. Thus the common velocity is constant on $[0,6]$, which implies that $\oa_i(\cdot)\equiv\oa_i$ on this interval. Hence $\eta_{12}(\cdot)\equiv\eta_{12}$ on $[0,6]$. It allows us to rewrite (2) as
\begin{equation*}
\dot\ox_1(t)=\left(3\sqrt2\oa_1-\frac{\sqrt2}{2}\eta_{12},-3\sqrt2\oa_1+\frac{\sqrt2}{2}\eta_{12}\right),\;
\dot\ox_2(t)=\left(\frac{3\sqrt2}{2}\oa_2+\frac{\sqrt2}{2}\eta_{12},-\frac{3\sqrt2}{2}\oa_2-\frac{\sqrt2}{2}\eta_{12}\right),
\end{equation*}
which gives ua by integration the trajectory formulas
\begin{equation*}
\left\{\begin{array}{ll}
\ox_1(t)=\left(-48-\frac{6}{\sqrt2}+\left(3\sqrt2\oa_1-\frac{\sqrt2}{2}\eta_{12}\right)t,48+\frac{6}{\sqrt2}+\left(-3\sqrt2\oa_1+\frac{\sqrt2}{2}\eta_{12}\right)t\right),\\
\ox_2(t)=\left(-48+\left(\frac{3\sqrt2}{2}\oa_2+\frac{\sqrt2}{2}\eta_{12}\right)t,48+\left(-\frac{3\sqrt2}{2}\oa_2-\frac{\sqrt2}{2}\eta_{12}\right)t\right).
\end{array}\right.
\end{equation*}
We have $\dot\ox_1(t)=\dot\ox_2(t)$ on $[0,6]$ due to the same velocity of the the participants, which yields
\begin{equation}\label{e:eta cal}
\eta_{12}=\dfrac{6\oa_1-3\oa_2}{2}.
\end{equation}
Furthermore, it follows from the optimality conditions in (5) and (8) that
\begin{equation}\label{e:a-q}
\lm\oa_1=-3\sqrt2q^x_{11}(t)+3\sqrt2q^x_{12}(t),\;\lm\oa_2=-\frac{3\sqrt2}{2}q^x_{21}(t)+\frac{3\sqrt2}{2}q^x_{22}(t)\;\mbox{ on }\;[0,6].
\end{equation}
If $\eta_{12}>0$, we immediately deduce from (4) that
\begin{equation*}
-q^x_{11}(t)+q^x_{12}(t)=-q^x_{21}(t)+q^x_{22}(t)\;\mbox{ on }\;[0,6].
\end{equation*}
Combining the latter with \eqref{e:a-q} shows that $\oa_1=2\oa_2$ provided that $\lm>0$; otherwise we do not have enough information to proceed. Consider now the following two cases:\\
{\bf Case~1:} $\eta_{12}=0$. Then \eqref{e:eta cal} tells us that $\oa_2=2\oa_1$ and the cost functional reads as
$$
J[\oa_1]=1311\oa_1^2-36(96\sqrt2+6)\oa_1+\left(48+\dfrac{6}{\sqrt2}\right)^2+48^2
$$
while attaining its minimum at $\oa_1=\dfrac{(96\sqrt2+6)18}{1311}\approx 1.95$. Thus $\oa_2=2\oa_1\approx 3.9$, the minimum cost is $J\approx 66.49$, and the optimal trajectories are
calculated by
\begin{eqnarray}\label{opt-tr}
\ox_1(t)=\left(-48-\dfrac{6}{\sqrt2}+8.27t,48+\dfrac{6}{\sqrt2}-8.27t\right),\;\ox_2(t)=\left(-48+8.27t,48-8.27t \right).
\end{eqnarray}
{\bf Case 2:} $\eta_{12}>0$. By $\oa_1=2\oa_2$ in this case we get $\eta_{12}=\frac{9}{2}\oa_2$. The cost functional reads as
$$
J[\oa_2]=2040\oa^2_2-45(96\sqrt2+6)\oa_2+\left(48+\frac{6}{\sqrt2}\right)^2+48^2
$$
with its minimum value achieved at $\oa_2=\dfrac{45(96\sqrt2+6)}{4080}\approx1.56$. Hence $\oa_1=2\oa_2\approx3.12$, the minimum cost is $J\approx45.9$, and the optimal trajectories are given by the same formulas \eqref{opt-tr} as in Case~1. Observe that in both cases we have \eqref{u-cont} for the corresponding $u$-control $\ou(\cdot)=(\ou_1(\cdot),\ou_2(\cdot))$ satisfying the constraints $r_1\le\|\ou(t)\|\le r_2$, and thus we obtain a complete solution of the problem under consideration.
\end{example}\vspace*{-0.3in}

\section{Concluding Remarks}
\setcounter{equation}{0}\vspace*{-0.1in}
In this paper we formulate a new class of optimal control problems for discontinuous differential inclusions governed by a nonconvex controlled version of the sweeping process with pointwise control and intrinsic state constraints. This class of problems is motivated by applications to optimal control of dynamical systems arising in the planar crowd motion model while being important for its own sake as well for other applications. The main thrust of the paper is developing the method of discrete approximations to derive necessary optimality conditions for the so-called ``intermediate" local minimizers in this highly nonstandard and challenging class of optimal control problems. This is done here by establishing well-posedness and strong convergence of discrete approximations together with rich calculus and explicit computations of robust first-order and second-order generalized differential constructions of variational analysis. In this way we obtain a rather comprehensive set of necessary optimality conditions expressed entirely via the problem data and then illustrate their applications by two examples related to real-life models. More complete applications to the controlled planar crowd model are presented in \cite{cm-crowd}. We also plan to apply the obtained results and their further developments to various models of the sweeping type arising in mechanics, engineering systems with hysteresis, robotics, etc.

\end{document}